\documentclass[11pt]{article}

\usepackage{amsfonts}
\usepackage{amssymb}
\newtheorem{st}      {Theorem}
\newtheorem{prop}  {Proposition}
\newtheorem{lem} {Lemma}
\newtheorem{cor} {Corollary}

\newtheorem{defin}     {Definition}
\newtheorem{exam}{Examples}
 
\newtheorem{rmk} {Remark}

\newcommand{\rmap}{\longrightarrow}
\newcommand{\Boxe}{\raisebox{.8ex}{\framebox}}
\newcommand{\lmap}{\longleftarrow}

\newcommand{\U}{\ensuremath{\mathcal{U}}}

\newcommand{\Ka}{\ensuremath{K}}
\newcommand{\Ha}{\ensuremath{H}}
\newcommand{\G}{\ensuremath{G}}
\newcommand{\xx}{\ensuremath{\mathcal{X}}}

\newcommand{\A}{\ensuremath{\mathcal{A}}}

\newcommand{\F}{\ensuremath{\mathcal{F}}}

\newcommand{\ps}{{\raise 1pt\hbox{\tiny (}}}

\newcommand{\pss}{{\raise 1pt\hbox{\tiny [}}}
\newcommand{\pdd}{{\raise 1pt\hbox{\tiny ]}}}
\newcommand{\pd}{{\raise 1pt\hbox{\tiny )}}}

\newcommand{\bs}{{\raise 1pt\hbox{\tiny [}}}
\newcommand{\bd}{{\raise 1pt\hbox{\tiny ]}}}
\newcommand{\nG}[1]{\ensuremath{\G^{\ps #1\pd}}}
\newcommand{\nK}[1]{\ensuremath{K^{\ps #1\pd}}}

\def\cross{\mathinner{\mathrel{\raise0.8pt\hbox{$\scriptstyle>$}}
                 \joinrel\mathrel\triangleleft}}
\def\dcross{\mathinner{\triangleright
                 \triangleleft}}
\newcommand{\nH}[1]{\ensuremath{H^{\ps #1\pd}}}

\def\compose{{\raise 1pt\hbox{$\scriptscriptstyle\circ$}}}

\advance\textheight by \topskip
\usepackage[all]{xy}
\usepackage{epsf}
\usepackage{amsfonts}
\usepackage{amssymb}
\usepackage{amsmath}

\begin{document}

\title{{\bf Differentiable and algebroid cohomology, Van Est isomorphisms,
and characteristic classes}\thanks{Research supported by NWO}\rm \\
$\text{}$}
\author {by Marius Crainic}
\pagestyle{myheadings}
\date{} 
\maketitle
\begin{abstract}
\hspace*{-.2in}In the first section we discuss Morita invariance of differentiable/algebroid cohomology.\\
In the second section we extend the Van Est isomorphism to groupoids.
As a first application we clarify the connection between differentiable 
and algebroid cohomology (proved in degree $1$, and conjectured in degree $2$ by Weinstein-Xu \cite{WeXu}).
As a second application we extend Van Est's argument for the integrability of Lie algebras. Applied to Poisson manifolds,
this immediately implies the  integrability criterion of Hector-Dazord  \cite{DaHe}.\\
In the third section we describe the relevant characteristic classes of representations, living in algebroid cohomology,
as well as their relation to the Van Est map. 
This extends Evens-Lu-Weinstein's characteristic class $\theta_{L}$ \cite{ELW} (hence, in particular, the modular class of Poisson manifolds),
and also the classical characteristic classes of flat vector bundles \cite{BiLo, KT}.\\
In the last section we describe applications to Poisson geometry.\\

\end{abstract}

\tableofcontents



\section*{Introduction}

\hspace*{.3in}The classical Van Est isomorphism \cite{vEst, VEst, EstCartan} is the main tool which relates
the differentiable cohomology of Lie groups to the cohomology of
Lie algebras. This paper grew out of author's attempts to
understand the differentiable cohomology of groupoids, the
cohomology of Lie algebroids, as well as the relations between
them.\\
\hspace*{.3in}Lie groupoids/algebroids \cite{McK} have proven to be very useful in several areas: foliation theory \cite{CM, Haefl, CoOp} Poisson geometry \cite{Vais, WeXu}, non-commutative geometry \cite{CoOp}, analysis on singular spaces \cite{Nistor, WeNi}, geometry of connections etc. 
Roughly speaking, a Lie groupoid consists of a base manifold (the space of
objects) and a ``group of symmetries'' (the space of arrows). Lie algebroids are
the infinitesimal version of Lie groupoids. In contrast with Lie algebras, there is no Lie
third theorem for Lie algebroids. Moreover, it seems we are quite far from a complete understanding of
this failure. Positive integrability results are however very important, e.g. for geometric quantization, the geometry of foliations, analysis on singular spaces.\\
\hspace*{.3in}Although differentiable/algebroid
cohomology are straightforward extensions of the corresponding
notions for Lie groups, they are extremely relevant when applied
to the list of examples above. Particular cases are: De Rham
cohomology, Lie algebra cohomology, Poisson cohomology, foliated cohomology, certain
cohomology groups of classifying spaces. Despite this, very little is known about the properties/relevance
of these cohomologies for general Lie groupoids. What basic properties do they enjoy (e.g. are they Morita invariant)? Which is the relation between differentiable and algebroid cohomology? Which 
are the invariants which live in these cohomologies? Apart from being basic questions on the theory of Lie groupoids/algebroids, they are also relevant to the applications of the general theory. For instance, in Poisson geometry, the relation between differentiable and algebroid cohomology is relevant
in the process of quantizing Poisson manifolds \cite{WeXu}. 
In the same direction, one knows \cite{ELW} that the modular class of
a Poisson manifolds lives in the world of Lie algebroids,
and appears as the (first) characteristic class of a one dimensional representation. 
In non-commutative geometry, the cyclic cohomology of convolution
algebras is undoubtedly  related to differentiable cohomology
of groupoids (this is clear for instance from \cite{Crath, Niii}); 
in this direction, note that one of the missing
steps for extending the results of \cite{Niii}  from Lie groups to
Lie groupoids is the lack of Van Est-type techniques for groupoids. Also, there are clear 
indications that cohomological methods might be useful to the integrability problem. The best example
is probably Van Est's proof \cite{VEst} of Lie's third fundamental theorem.\\
\hspace*{.3in}The purpose of this paper is to study these general properties of differentiable/algebroid cohomology,
and to describe some applications. Below we give an outline of the 
main results, as well as of the relation with other known results/conjectures.  Note that all the applications to Poisson geometry have been collected in the last section. \\

\hspace*{.1in}{\it Morita invariance:}  Intuitively, two
groupoids are Morita equivalent if they have the same space of
orbits (i.e. the same transversal geometry). This notion is important in the theory
of foliations (where groupoids are often replaced by Morita equivalent ones; see e.g. \cite{CM, Haefl}), non-commutative geometry (recall that the $C^*$-algebras defined by two Morita equivalent groupoids are stably isomorphic; see  \cite{CoOp, MRW}), Poisson geometry (recall that Morita equivalent Poisson manifolds can be integrated by Morita equivalent symplectic groupoids; see \cite{XuM}), and to the geometry of principal bundles (gauge groupoids are Morita equivalent to Lie groups). First we prove Morita invariance of differentiable cohomology (Theorem \ref{Moritatheorem}).
Since the notion of Morita equivalence of Lie algebroids is a bit problematic (there are several natural, but non-equivalent definitions), we restrict ourselves to an invariance of algebroid cohomology which we find relevant to applications.
This is our second theorem (Theorem \ref{fibers}).\\

\hspace*{.1in}{\it Van Est isomorphisms:} The next theorem is an extension of Van Est's isomorphism to groupoids
(Theorem \ref{VEth}; see also the comments at the beginning of Section \ref{AVEi}). This immediately implies a version 
of Haefliger's conjecture \cite{difcoh, conj} for differentiable cohomology (Corollary \ref{Haeflconj}). As a first application of Theorem \ref{VEth} we clarify the relation between differentiable/algebroid cohomology, via a Van Est map
\begin{equation}\label{intr2}  
\Phi: H^{*}_{d}(\G)\rmap H^{*}(\mathfrak{g}) \ .
\end{equation}
This is the object of our Theorem \ref{theWeXu}. In degree one it becomes a (slight improvement of a) theorem of Weinstein-Xu, and in degree two it proves a conjecture by the same authors (see Theorem  $1.3$ and the
comments at page $172$ of \cite{WeXu}). In combination with the Morita invariance, this also clarifies the
invariance of Poisson cohomology (see Corollary \ref{abcMP}), which has been known in degree one only \cite{GiLu}.  \\

\hspace*{.1in}{\it An integrability result:} In Theorem \ref{integration} we present another application of Van Est's isomorphism. We prove that, given an extension of Lie algebroids
\[ 0\rmap \mathfrak{l}\stackrel{i}{\rmap} \mathfrak{h} \stackrel{\pi}{\rmap}\mathfrak{g} \rmap 0 \]
such that $\mathfrak{l}$ is abelian, and $\mathfrak{g}$ admits an integration which is $\alpha$-two-connected, then $\mathfrak{h}$ is integrable. Note that, in the case of Lie algebras, this immediately implies Lie's third
theorem (and this is Van Est's proof \cite{VEst} mentioned above). Indeed, if $\mathfrak{h}$ is an arbitrary Lie algebra, we can take $\mathfrak{l}$ to be the center of $\mathfrak{h}$, and, since $\mathfrak{g}\subset gl(\mathfrak{g})$ by the adjoint 
representation, $\mathfrak{g}$ is integrable. This shows that, quite surprisingly, the non-trivial part of Van Est's argument does extend to groupoids 
(see also our Remarks \ref{commint}). There is also strong evidence that these Van Est-techniques, combined with Cattaneo-Felder approach
to integrability \cite{Cat} could clarify the integrability problem (compare to the proof of Lie's third in \cite{Duis}). \\   
\hspace*{.3in}As an immediate consequence of our integrability result we obtain a slight improvement and a more conceptual proof  of Dazord-Hector's criterion \cite{DaHe}
for the integrability of regular Poisson manifolds (see our Corollary \ref{corHD}). This shows that Dazord-Hector's result is precisely
Van Est's argument applied to Poisson manifolds. \\
\hspace*{.3in}Comparing with another known integrability result, the one of Moerdijk-Mrcun \cite{MM}
and Nistor \cite{Nistor}, 
our Theorem \ref{integration} is both stronger (since we allow Lie algebroids
more general then semi-direct products), and weaker (since we have
to assume abelian algebroids instead of bundles of Lie algebras). An improvement which contains both 
integrability results would probably give a much better understanding on Lie's third theorem for groupoids.\\

\hspace*{.1in}{\it Characteristic classes in algebroid cohomology:} The last problem that we take concerns the
invariants living in algebroid cohomology. The conclusion
is that, given a $n$-dimensional representation $E$ of a Lie
algebroid $\mathfrak{g}$, the Chern classes of $E$ viewed in
$H^*(\mathfrak{g})$ must vanish (Theorem \ref{vanishing}), while new
classes $u_{2k-1}(E)\in H^{2k-1}(\mathfrak{g})$, $1\leq k\leq n$ 
show up. Note the analogy with Bott's vanishing theorem, and
the construction of secondary characteristic classes for foliations \cite{Bo}. We obtain a characteristic map 
\begin{equation}\label{intr1} 
U: Rep(\mathfrak{g})\rmap \mathbb{Z}\ltimes H^{odd}(\mathfrak{g}) \ . 
\end{equation}
See also \cite{Ginz}. Part of our motivation was to extend the construction
of the characteristic class $\theta_{L}$ found by Evens-Lu-Weinstein \cite{ELW}
(hence, in particular, of the modular class of a Poisson manifold), which we recover
when $n= k= 1$. Although we also sketch an explicit approach (in the spirit of \cite{BiLo}), the way we introduce the higher 
$u_{2k-1}$'s is by extending the well-known construction \cite{KT} of 
characteristic classes of flat vector
bundles over manifolds (which we recover in the case of tangent
Lie algebroids).  One
advantage of our (Chern-Weil type) approach is that it can be used in the presence of other structural groups (not
only $GL_n$). Since some possible applications of this
construction are to Lie algebroids which are not integrable, or
whose integrations are difficult to work with, we do everything at
the algebroid level. However, when $E$ is the representation of a
Lie groupoid, we show (Theorem \ref{RwVE}) that its characteristic
classes (\ref{intr1}) come from the differentiable cohomology via the Van Est map (\ref{intr2}).
An immediate consequence of this is the Morita invariance of the modular class of 
Poisson manifolds (Corollary \ref{corolar3}),
which has been known under certain conditions only (see Theorem $4.2$ and the conjecture $4.6$ in \cite{Gi}).\\

\hspace*{.1in} {\bf Acknowledgments}: All the discussion I had with K. Mackenzie and
I. Moerdijk were an essential source of inspiration. Useful were also the comments that 
V. Ginzburg and R.L. Fernandes had on a preliminary version of this paper. Regarding R.L. Fernandes 
I have learned, in a late stage of my work, about his construction of characteristic classes for algebroids \cite{RLF}
(and for Poisson manifolds in particular \cite{Fer2}).
Recently we have explained the connection between our approaches \cite{Crai1, Crai2} (see also our Examples \ref{fernandes}). 
I have also benefited from the discussions with J. Mrcun. 




\section{Differentiable and algebroid cohomology}

\hspace*{.3in}This section is an exposition of basic definitions
and properties concerning groupoids, differentiable cohomology,
Lie algebroids, and algebroid cohomology \cite{difcoh, Haefl, McK, McHi}. It is a combinations of well-known
definitions which we recall for reference, and some remarks which we find
important for understanding the objects under discussion. Here we prove 
Morita invariance of differentiable cohomology (Theorem \ref{Moritatheorem}), and an
invariance for algebroid cohomology (Theorem \ref{fibers}).


\subsection{Groupoids}
\label{basicgroupoids}


\hspace*{.3in}Recall that a {\it groupoid} $\G$ is a (small) category in which every arrow is
invertible. We will write $\nG{0}$ and $\nG{1}$ for the set of
objects and the set of arrows in $\G$, respectively. We also say that $\G$ is a groupoid
over $\nG{0}$. 
The source and target maps are denoted by $\alpha, \beta: \nG{1}\rmap \nG{0}$,
while $m(g, h)= g\compose h$ is the composition, and $i(g)= g^{-1}$ denotes the inverse of $g$.
One calls $\G$ a  {\it Lie groupoid} if $\nG{0}$ and
$\nG{1}$ are smooth manifolds, all the structure maps are smooth, and $\alpha$ and $\beta$ are submersions.
We denote by $\nG{p}$ the space of $p$-composable strings 
\begin{equation}\label{string} 
x \stackrel{g_1}{\lmap} \stackrel{g_2}{\lmap} \ldots \stackrel{g_p}{\lmap}
\end{equation}
and by
\begin{equation}\label{betamap} 
\beta= \beta_p: \nG{p}\rmap \nG{0} 
\end{equation}
the map which associates to a string (\ref{string}) the element $x$.\\
\hspace*{.3in}A {\it left action} of $\G$ on a manifold $X$ consists of a smooth map $\epsilon: X\rmap \nG{0}$
called {\it the moment map} of the action, together with a smooth map $\nG{1}* X\rmap X$, $(g, x)\mapsto gx$
defined on the space of pairs $\{(g, x): \alpha(g)= \epsilon(x)\}$, which satisfies the usual identities
of an action. One defines the action groupoid $\G\cross X$ as the groupoid over $X$ with $\nG{1}* X$ 
as space of arrows, the multiplication as source map, the second projection
as target map, and multiplication $(g, x)(h, y)$ $=$ $(gh, y)$ (defined when $x= hy$). 
A {\it left $\G$-bundle} consists of a (left) $\G$-space $P$, and a $\G$-invariant surjective submersion
$\pi: P\rmap B$. It is called {\it principal} if
\[ \G* P\rmap P\times_{B}P,\ \ (g, p)\mapsto (gp, p) \]
is a diffeomorphism. Similarly one defines the notion of a right action of $\G$ on a space $X$, the action
groupoid $X\cross \G$, and the notion of right $\G$-bundle.\\
\hspace*{.3in}A {\it morphism} between two Lie groupoids is a smooth functor. If one is interested in
the orbit space rather then on the groupoid itself (i.e. in ``transversal structures''), one has to
relax the notion of morphism. Recall \cite{Haefl, Mrcun} that a {\it generalized morphism}
$\phi: \G\rmap \Ha$ between two Lie groupoids is determined by a manifold $P$
(view it as ``the graph'' of $\phi$) endowed with a left $\G$-action with moment map denoted $\epsilon: P\rmap \nG{0}$,
a right $\Ha$-action, with moment map denoted $\eta: P\rmap \nH{0}$, such that $\epsilon: P\rmap \nG{0}$ is a principal $\Ha$-bundle,
$\eta: P\rmap \nH{0}$ is a $\G$-bundle, and the two actions are compatible. 
The set of isomorphism classes of such $P$'s forms precisely the set $Hom_{gen}(\G, \Ha)$
of generalized homomorphisms from $\G$ to $\Ha$. We say that $P$ defines 
a {\it Morita equivalence} if $P$ is principal also as a $\G$-bundle. Given $\phi: \G\rmap \Ha$ represented by the bundle $P$, $\psi: \Ha\rmap \Ka$ represented by
the bundle $Q$, one can define the
composition $\psi\compose \phi$ represented by the bundle $P\otimes_{\Ha}Q$, which
is the quotient of $P\times_{\nH{0}}Q$ by the diagonal action
of $\Ha$. One can easily see that Morita equivalences are precisely
the isomorphisms of the resulting category. \\
\hspace*{.3in}Note that any morphism $\phi: \G\rmap \Ha$ can be viewed as a generalized morphism 
whose graph $P_{\phi}$ consists of pairs $(x, h)\in \nG{0}\times \nH{1}$ such that $\phi(x)= \beta(h)$,
with $\epsilon(x, h)= x$, $\eta(x, x)= \alpha(h)$, and the obvious actions. We say that $\phi$ is an {\it essential equivalence}
if $P_{\phi}$ defines a Morita equivalence.\\

\begin{exam}\label{exam1}\emph{ Any Lie group is a groupoid with only one object,
and any manifold can be viewed as a groupoid with identity arrows only.
Any surjective submersion
$\pi: M\rmap B$ induces a groupoid $M\times_BM$ over $M$, consisting
on pairs $(x, y)\in M\times M$ so that $\pi(x)= \pi(y)$, with the projections as
source and target, and with the obvious multiplication.
One has an obvious Morita equivalence $M\times_{B}M\cong B$ with
$M$ as defining bundle. The particular case where $B$ is a point
gives the pair groupoid $M\times M$. The fundamental groupoid of a manifold $M$
is Morita equivalent to the fundamental group of the manifold (take the
universal cover of $M$ as defining bundle).  Any transitive groupoid is Morita equivalent
to a Lie group. All these can be viewed as an instance of the gauge groupoid: given any groupoid $\G$, and any principal $\G$-bundle $P\rmap B$,
one can form the gauge groupoid over $P$, which is $P\otimes_{\G}P$ (we use the notations above), with the obvious structure maps;
$P$ gives a Morita equivalence $P\otimes_{\G}P\cong \G$.
In (transversal) foliation theory, the holonomy and the monodromy groupoids,
as well as Haefliger's groupoids play a central role. The former one is etale (i.e.
its source is a local diffeomorphism), while the first two are Morita equivalent to etale ones.
In general, any foliation groupoid is Morita equivalent to an etale one \cite{CM}.\\
\hspace*{.3in}Although important examples coming from foliation theory may be non-Hausdorff,
in order to deal effectively with the differentiable cohomology of groupoids, we need
to make the following}
\end{exam}

{\bf Overall assumption:} All groupoids $\G$ in this paper are assumed to be Hausdorff,
in the sense that $\nG{0}$ and $\nG{1}$ are Hausdorff manifolds.


\subsection{Differentiable cohomology and Morita invariance}
\label{DcMi}

\hspace*{.3in}Let $\G$ be a Lie groupoid. A {\it left representation} of $\G$ is a (real or complex) vector
bundle $\pi: E\rmap \nG{0}$, endowed with an action of $\G$ whose moment map is precisely
$\pi$, and which is fiberwise linear. Denote by $Rep(\G)$ the semi-ring of isomorphism classes 
of (left) representations of $\G$. An $E$-valued differentiable $p$-cochain on $\G$ 
is a smooth map $c$ which associates to a string (\ref{string})
an element $c(g_1, \ldots, g_p)\in E_x$. The space $C^{p}_{d}(\G; E)$
of $p$-cochains coincides with the space of smooth sections of $\beta_{p}^{*}E$,
and comes
equipped with a differential $d: C^{p}_{d}(\G; E)\rmap C^{p+1}_{d}(\G; E)$,
\begin{eqnarray} 
(dc)(g_{1}, \ldots , g_{p}, g_{p+1}) & = & g_1 c(g_{2}, \ldots , g_{p+1}) + \\
   & + & \sum_{i=1}^{p} (-1)^i c(g_{1}, \ldots , g_{i}g_{i+1}, \ldots , g_{p+1})+ (-1)^{p+1}c(g_{1}, . . . , g_{p}) .
\end{eqnarray}
The {\it differentiable cohomology} of $\G$ with coefficients in $E$, denoted $H^{*}_{d}(\G; E)$,
is defined as the cohomology of the resulting complex $(C^{*}_{d}(\G; E), d)$. When $E$ is the trivial line bundle, we simplify the notations to $C^{*}_{d}(\G)$ and $H^{*}_{d}(\G)$.
The usual cup-product:
\[ (c_1\cup c_2)(g_1, \ldots g_{p+q})= c_1(g_1, \ldots g_p) \otimes (g_1\ldots g_p)c_2(g_{p+1}, \ldots, g_{p+q}) ,\]
defines a product structure $C_{d}(\G; E)\otimes C_{d}(\G; F)$ $\rmap$ $C_{d}(\G; E\otimes F)$ which passes to
cohomology because it satisfies the Leibniz identity
$d(c_1\cup c_2)= d(c_1)\cup c_2 + (-1)^{deg(c_1)} c_1\cup d(c_2)$.
In particular,  $(C^{*}_{d}(\G), d)$ is a DG algebra, $H^{*}_{d}(\G)$
is a (graded) algebra, and $H^{*}_{d}(\G; E)$ are $H^{*}_{d}(\G)$-modules. \\
\hspace*{.3in}The notion of representation is a basic example of transversal structure,
in the sense that any Morita equivalence $\phi: \G\rmap \Ha$ induces an isomorphism
$\phi^*: Rep(\Ha)\tilde{\rmap} Rep(\G)$. More generally, any generalized morphism 
$\phi\in Hom_{gen}(\G, \Ha)$ induces a map $\phi^*: Rep(\Ha)\rmap Rep(\G)$, natural
on $\phi$. If $E\in Rep(\Ha)$, and $\phi$ is represented by the bundle $P$, then the pull-back 
$\eta^*E$ of the vector bundle $E$ to $P$ is equipped with a free right $\Ha$-action, and
a tautological left action of $\G$. Define then $\phi^*E:= \eta^*(E)/\Ha$ which is a vector bundle 
over $P/\Ha= \nG{0}$, equipped with a left action of $\G$, i.e. a representation of $\G$. 

\begin{st}\label{Moritatheorem} Any Morita equivalence $\phi: \G\rmap \Ha$ induces isomorphisms:
\[ \phi^*: H^{*}_{d}(\Ha)\tilde{\rmap}H^{*}_{d}(\G) \ .\]
\hspace*{.3in}More generally, for any $\phi\in Hom_{gen}(\G, \Ha)$, and
any $E\in Rep(\Ha)$, there is an induced homomorphism
\[ \phi^*: H^{*}_{d}(\Ha; E)\rmap H^{*}_{d}(\G; \phi^{*}E) \ .\]
The construction is natural on $\phi$, and is compatible with the product structure.
If $\phi$ is a morphism of groupoids, then $\phi^*$ is just the map induced by the composition with $\phi$.
\end{st}

{\it Proof:} We prove the theorem with trivial coefficients
(in the general case, there are obvious modifications). 
Let $P$ be a bundle representing $\phi$. 
For each $p, q$, $P\times_{\nH{0}}\nH{q}$ is a (left) $\G$-space,
$\nG{p}\times_{\nG{0}}P$ is a (right) $\Ha$-space. Remark that the spaces of composable arrows 
of the associated crossed product groupoids are related by
\[ (\G\cross (P\times_{\nH{0}}\nH{q}))^{\ps p\pd}= ((\nG{p}\times_{\nG{0}}P)\cross \Ha)^{\ps q\pd}= 
\nG{p}\times_{\nG{0}}P\times_{\nH{0}}\nH{q} \ .\]
We form a double complex $C(P)$ with $C^{p, q}(P)= C^{\infty}(\nG{p}\times_{\nG{0}}P\times_{\nH{0}}\nH{q})$,
and with the differentials defined so that the $q^{th}$ row $C^{*, q}(P)$ is the
complex computing the differential cohomology of $\G\cross (P\times_{\nH{0}}\times\nH{q})$,
and similarly, the $p^{th}$ column $C^{p, *}(P)$ is the complex computing
the differentiable cohomology of $(\nG{p}\times_{\nG{0}}P)\cross \Ha$. 
The complexes computing the differentiable cohomologies of $\G$ and $\Ha$
come as co-augmentations of the columns, and rows, respectively:
\[ C_{d}(\G) \stackrel{\epsilon^*}{\rmap} C(P) \stackrel{\eta^*}{\lmap} C_{d}(\Ha) \ .\]
We show that the condition that $P$ is principal as an $\Ha$ space
implies that the co-augmented columns of $C(P)$ are acyclic, hence $\epsilon^*$ 
induces isomorphisms in cohomology. This will give the desired map in cohomology,
\[ \phi^*: H_{d}^{*}(\G) \stackrel{(\epsilon^*)^{-1}}{\rmap} H^*(C(P)) \stackrel{\eta^*}{\rmap} H^{*}_{d}(\Ha) ,\]
and, if $P$ is principal also as a $\G$-space (i.e. is a Morita equivalence),
a similar argument proves that also $\eta^*$ induces isomorphisms in cohomology. 
Since each of the groupoids $(\nG{p}\times_{\nG{0}}P)\cross \Ha$ is diffeomorphic 
to a groupoid of type $X\times_BX$ (choose $X= \nG{p}\times_{\nG{0}}P$, $B= \nG{p}$), the acyclicity of the columns
follows from the following:

\begin{lem}For any surjective submersion $\pi: X\rmap B$,
and any vector bundle $E$ on $B$,  
the differentiable cohomology $H^{*}_{d}(X\times_BX; \pi^*E)$ is zero
in positive degrees, and is $C^{\infty}(B; E)$ in degree zero.
\end{lem}

{\it Proof:} We have to prove that the complex $C^*= C_{d}^{*}(X\times_{B}X)$, together with the co-augmentation $\pi^*: C^{\infty}(B)$ $\rmap$ $C^{\infty}(X)= C^0$, is exact. If $\pi$ admits a continuous section $s: B\rmap X$, we have the explicit homotopy 
$h(c)(x_1, \ldots , x_{n-1})= c( s(pi(x_1)), x_2, \ldots , x_{n-1})$. In the general case, we choose an open
cover $\U$ of $B$ over which $\pi$ admits local sections, and we use a Mayer-Vietoris argument \cite{BoTu}. \ $\Boxe$\\

\hspace*{.1in}To prove that $\phi^*$ is compatible with the product structure, it suffices
to find a (bigraded) product on $C(P)$, such that $\epsilon^*$ and $\eta^*$
are compatible with the products. For $c\in C^{p, q}$, $c\,'\in C^{p\,', q\,'}$ we define
$\omega\cdot\eta$ $\in$ $C^{p+q, p\,'+ q\,'}$ by
\begin{eqnarray}
 & (c\cdot c\,')(g_1, \ldots , g_{p+ p\,'}, x, h_{1}, \ldots , h_{q+ q\,'}) = \nonumber\\
 & c(g_1, \ldots , g_p, g_{p+1}\ldots g_{p+p\,'}x, h_1, \ldots , h_q) 
   c\,'(g_{p+1}, \ldots , g_{p+p\,'}, xh_1\ldots h_{q}, h_{q+1}, \ldots , h_{q+q\,'}) .\nonumber\
\end{eqnarray}
\hspace*{.3in}If $\phi$ is a morphism of groupoids, denote by $\tilde{\phi}: C_{d}(\Ha)\rmap C_{d}(\G)$ the composition by 
$\phi$. To prove that $\tilde{\phi}$ induces $\phi^*$ in cohomology, it suffices to find a chain map
$\Phi: C(P)\rmap C_{d}(\G)$ (where $P= P_{\phi}$, see subsection \ref{basicgroupoids}) such that
$\Phi\compose \epsilon^*= Id$, and $\Phi\compose \eta^{*}= \tilde{\phi}$. For $c\in C^{p, q}(P_{\phi})$ $=$ $C^{\infty}(\nG{p}\times_{\nH{0}}\nH{q+1})$ we set
\[ \Phi(c) (g_1, \ldots , g_{p+q})= c(g_1, \ldots, g_p, 1_{\phi\ps \alpha\ps g_p\pd\pd}, \phi(g_{p+1}) ,\ldots , \phi(g_{p+q})) .\]
\hspace*{.3in}We now prove the naturality $\phi^*\psi^*= (\psi\phi)^*$ for 
$\phi\in Hom_{gen}(\G, \Ha)$, $\psi\in Hom_{gen}(\Ha, \Ka)$.
Choose $P$ and $Q$ representing $\phi$, and $\psi$, respectively, 
and we consider a triple complex $C(P, Q)$
with
\[ C^{p, q, r}(P, Q)= C^{\infty}(\nG{p}\times_{\nG{0}}P\times_{\nH{0}}\times\nH{q}\times_{\nH{0}}Q\times_{\nK{0}}\nK{r}) \]
defined analogous to $C(P)$. The double complexes $C(P)$, $C(Q\compose P)$, and $C(Q)$ appear all
as co-augmentations of $C(P, Q)$, and the first two of this co-augmentations are quasi-isomorphisms. 
Moreover, they are compatible with the co-augmentations of $C(P)$, $C(Q\compose P)$, and $C(Q)$,
and then the conclusion follows by a diagram chasing.\ $\Boxe$\\

\begin{exam}\label{example1}\emph{ When $G$ is a Lie group, one recovers the ususal differentiable cohomology
of Lie groups. In particular, for discrete groups, one obtains the usual group cohomology.
By Morita invariance, one can compute the differentiable cohomology in many examples:
transitive groupoids, fundamental groupoids, foliation groupoids. Related to the last class
of examples, remark that if $\G$ is a Hausdorff etale groupoid, then $H^{*}_{d}(\G)$
is (by definition) the same as Haefliger's cohomology \cite{difcoh} with coefficients in the
sheaf $\A$ of smooth functions. Hence one obtains a sheaf $\tilde{A}$
over the classifying space $B\G$ \cite{Se} so that $H^{*}_{d}(\G)\cong H^{*}(B\G; \tilde{A})$
(this was conjectured by Haefliger and proved in \cite{conj}).}
\end{exam}

\subsection{From Lie groupoids to Lie algebroids}
\label{FGA}

\hspace*{.3in}Analogous to the construction of the Lie algebra of
a Lie group, one can construct the infinitesimal version of a Lie groupoid $\G$.
We recall here some of these constructions, and we refer to \cite{McK} for more details. 
The central role is played by the vector bundle $\mathfrak{g}$ which 
is defined as the restriction
along $u: \nG{0} \hookrightarrow \nG{1}$ of the vector bundle $T^{\alpha}(\nG{1}) =
Ker(d\alpha: T\nG{1} \rightarrow T\nG{0})$ of ``$\alpha$-vertical'' tangent vectors on $\nG{1}$.
The fiber of $\mathfrak{g}$ at $x\in \nG{0}$ is the tangent space
at $1_x$ of $G(x, -) = \alpha^{-1}(x)$. Moreover, the differential of the target $\beta$ of $\G$
induces a map of vector bundles over $\nG{0}$, {\it the anchor map}, 
\[ \rho: \mathfrak{g}\rmap TM .\]
Any section $X\in \Gamma(\mathfrak{g})$ defines
a vector field $\tilde{X}$ on $\nG{1}$ by the formula 
\[ \tilde{X}(g)= (dR_{g})_{x}( X(x))\ ,\]
where $R_{g}: \G(x, -)\rmap \G(y, -)\subset \nG{1}$ is the right multiplication by $g$,
$x= \beta(g)$, $y= \alpha(g)$. Denoting by $\xx^{\alpha}(\G)= \Gamma(T^{\alpha}(\nG{1})$ 
the space of ``$\alpha$-vertical'' vector fields on $\nG{1}$, and by 
$\xx^{\alpha}_{inv}(\G)$ the subspace of right invariant
vector fields, i.e. vector fields $X\in \xx^{\alpha}(\G)$ with the property that
$X(gh)= (dR_h)_{g}(X(g))$ for all composable arrows $g, h$ of $\G$, the
construction above defines an isomorphism:
\[ \Gamma(\mathfrak{g}) \tilde{\rmap} \xx^{\alpha}_{inv}(\G)\subset \xx(\nG{1}), \ \ X\mapsto \tilde{X} \ .\]
Expressing the Lie brackets in terms of flows (see below) we see that $\Gamma(\mathfrak{g})\cong \xx^{\alpha}_{inv}(\G)$ is closed under the usual Lie bracket of vector fields on $\nG{1}$. Hence $\Gamma(\mathfrak{g})$ comes equipped with a Lie bracket $[\cdot \, , \cdot ]$, and it is not difficult to see that it is related to
the anchor map by the formula:
\begin{equation}\label{relalg} 
[X, fY] = f[X,Y]+ \rho(X)(f) \cdot Y 
\end{equation}
for all $X, Y \in \Gamma(\mathfrak{g})$ and $f \in C^{\infty}(\nG{0})$.
The resulting structure:
\[     (\mathfrak{g}, \rho, [\cdot \, , \cdot ])\]
is called the Lie algebroid of $\G$.\\
\hspace*{.3in}The connection between (the sections of) $\mathfrak{g}$ and $\G$ is given
by the flows. Given $X\in \Gamma(\mathfrak{g})$, we denote by $\phi_{X}$ the flow of
$\tilde{X}$ on $\nG{1}$, and, for $x\in \nG{0}$, put $\beta_{X}(t, x)= \beta\phi_{X}(t, x)$.
Then $\beta_{X}$ is a flow on $\nG{0}$, namely the flow of the vector field $\rho(X)$, while
the properties of $\tilde{X}$  show that 
\[ \phi_{X}(t, g)= \phi_{X}(t, x) g,\ \ \phi_{X}(t, x)\in \G(x, \beta_{X}(t, x)) \]
where $x= \beta(g)$. We then see that the Lie bracket in $\Gamma(\mathfrak{g})$ is given by
\begin{equation}\label{Liebracket}
[X, Y](x) = \frac{d}{dt}\frac{d}{ds}|_{t= s= 0} \phi_{X}(t, \beta_Y(s, \beta_X(-t, x)))\phi_Y(s, \beta_X(-t, x)) \phi_X(-t, x) .
\end{equation}
An action of $\G$ on a vector bundle $E$ over $\nG{0}$ has an infinitesimal version:
\begin{equation}\label{actiong} 
L_{X}(s)(x)= (\frac{d}{dt})_{ _{t=0}} \phi_X(t, x)^{-1} s(\beta_X(t, x)) \in E_x\ ,
\end{equation}
which defines a pairing $\Gamma(\mathfrak{g})\times \Gamma (E)\rmap \Gamma (E)$, $(X, s)\mapsto L_X(s)$.
These derivatives of sections of $E$ along sections of $\mathfrak{g}$ satisfy the basic relations: 
\begin{eqnarray}
\hspace*{.3in}  L_{fX}(s) & = & f L_{X}(s) \label{unu} \\
\hspace*{.3in}  L_X(fs) & = & fL_X(s)+ L_{X}(f)s \label{doi} \\
\hspace*{.3in}  L_{[X, Y]} & = & [L_{X}, L_{Y}] \label{trei}
\end{eqnarray}
for all $X, Y\in \Gamma(\mathfrak{g})$, $f\in C^{\infty}(\nG{0})$, $s\in \Gamma(E)$.

\subsection{Algebroids and their cohomology}
\label{ss1.3}

\hspace*{.3in}A {\it Lie algebroid} over a manifold $M$ is a triple
\[ (\mathfrak{g}, \, [\cdot \, , \cdot ] \, , \rho) \]
consisting of a vector bundle $\mathfrak{g}$ over $M$, a Lie bracket 
$[ \cdot \, , \cdot ]$ on the space $\Gamma(\mathfrak{g})$, and a
morphism of vector bundles $\rho: \mathfrak{g}\rmap TM$ ({\it the anchor}
of $\mathfrak{g}$), so that (\ref{relalg}) holds true for all
$X, Y \in \Gamma(\mathfrak{g})$ and $f \in C^{\infty}(M)$. \\
\hspace*{.3in}A {\it representation of $\mathfrak{g}$} is a vector bundle 
$\pi: E\rmap M$, together with a bilinear map (called the infinitesimal action
of $\mathfrak{g}$ on $E$)
\[ \Gamma(\mathfrak{g})\times \Gamma(E)\rmap \Gamma(E),\  (X, s)\mapsto L_{X}(s) ,  \]
satisfying the relations (\ref{unu})-(\ref{trei}) for all $X, Y\in \Gamma(\mathfrak{g})$, 
$f\in C^{\infty}(M)$, $s\in \Gamma(E)$. Isomorphism classes of representations of $\mathfrak{g}$, with
the direct sum, and tensor product as operations, form a semi-ring $Rep(\mathfrak{g})$. It actually comes endowed
with an involution $*$: $E^{*}$ is the dual of the conjugate of $E$ (hence, in the real case it is just the dual of $E$), with
$L_{X}(\omega)(s)= L_{X}(\omega(s))- \omega(L_X(s))$, for $\omega\in \Gamma(E^*)$, $s\in \Gamma(E)$. With this, a
metric on $E$ is called invariant if the induced isomorphism of vector bundles
$E^{*}\cong E$ is compatible with the action of $\mathfrak{g}$.\\
\hspace*{.3in}Of course, the motivating example is the Lie algebroid of
a Lie groupoid $\G$, and the representations induced by the representations
of $\G$.  As terminology, we call
integration of $\mathfrak{g}$, any Lie groupoid $\G$ so that $\mathfrak{g}= Lie(\G)$.
In contrast with the theory of Lie groups, not all Lie algebroids are
integrable. See however \cite{MM} for a large class of positive results. For instance one knows
that, if $\mathfrak{g}$ is integrable, then it has a unique integration $\G$ which is $\alpha$-simply connected
(i.e. whose $\alpha$-fibers are simply connected). Moreover, in this case any representation $E$ of $\mathfrak{g}$ 
can be uniquely integrated to a representation of $\G$ (this is e.g. a very particular case of Theorem $3.6$ in \cite{MM});
hence one has an isomorphism $Rep(G)\cong Rep(\mathfrak{g})$. \\

\hspace*{.1in}Let $E$ be a representation of the Lie algebroid $\mathfrak{g}$. An $E$-valued $p$-cochain on $\mathfrak{g}$ is a 
$C^{\infty}(M)$-multilinear antisymmetric map 
\[ \Gamma(\mathfrak{g})\times \ldots \times \Gamma(\mathfrak{g}) \ni (X_1, \ldots , X_p) \mapsto
\omega(X_1, \ldots , X_p)\in \Gamma(E) .\]
The space $C^{p}(\mathfrak{g}; E)$ of such cochains coincides
with the space of sections of the vector bundle $\Lambda^p\mathfrak{g}^*\otimes E$
over $M$, and comes equipped with a differential $d:C^{p}(\mathfrak{g}; E)\rmap C^{p+1}(\mathfrak{g}; E)$,
\begin{eqnarray}\label{differential}
d(\omega)(X_1, \ldots , X_{p+1}) & = & \sum_{i<j}
(-1)^{i+j-1}\omega([X_i, X_j], X_1, \ldots , \hat{X_i}, \ldots ,
\hat{X_j}, \ldots X_{p+1})) \nonumber \\
 & + & \sum_{i=1}^{p+1}(-1)^{i}
L_{X_i}(\omega(X_1, \ldots, \hat{X_i}, \ldots , X_{p+1})) .
\end{eqnarray}
Define {\it the cohomology of $\mathfrak{g}$ with coefficients in $E$}, denoted $H^{*}(\mathfrak{g}; E)$, as the
cohomology of the resulting complex $C^{*}(\mathfrak{g}; E)$. When $E$ is the
trivial line bundle (with the action $L_X(f)= \rho(X)(f)$, $X\in \Gamma(\mathfrak{g})$, $f\in C^{\infty}(M)$),
we simplify the notations to $C^{*}(\mathfrak{g})$ and $H^*(\mathfrak{g})$.\\
\hspace*{.3in}Note that, as in the case of differentiable cohomology, 
the usual wedge product (defined fiberwise) defines product structures 
$H^{*}(\mathfrak{g}; E)\otimes H^{*}(\mathfrak{g}; F)$ $\rmap$ $H^{*}(\mathfrak{g}; E\otimes F)$.
In particular, $C^{*}(\mathfrak{g})$ is a DG algebra, $H^*(\mathfrak{g})$ is a (graded) algebra,
and $H^*(\mathfrak{g}; E)$ are $H^*(\mathfrak{g})$-modules. Moreover, the usual Cartan calculus 
extends to $C^{*}(\mathfrak{g})$. More precisely, any $X\in \Gamma(\mathfrak{g})$ induces
Lie derivatives and interior products
\[ L_{X}: C^*(\mathfrak{g}) \rmap C^{*}(\mathfrak{g}), \ i_X: C^{*}(\mathfrak{g})\rmap C^{*-1}(\mathfrak{g})\ ,\]
\begin{eqnarray}
L_X(\omega)(X_1, \ldots X_p) & = & \sum \omega(X_1, \ldots, X_{i-1}, [X, X_i], X_{i+1}, \ldots , X_p) -\nonumber \\
                             & - & L_{X}(\omega(X_1, \ldots, X_p)) \nonumber \\
i_X(\omega)(X_1, \ldots , X_{p-1}) & = & \omega(X, X_1, \ldots , X_{p-1}) .\nonumber
\end{eqnarray}
The Lie derivatives are derivations of degree $0$, the interior products are derivations of degree $-1$,
and they satisfy the Cartan relations
\begin{eqnarray}
 \label{C1} di_{X}+ i_{X}d & = & L_{X}  \\       
 \label{C2}[L_{X}, L_{Y}] & = & L_{[X, Y]} \\                     
 \label{C3}[L_{X}, i_{Y}] & = & i_{[X, Y]} \\                
 \label{C4}[i_{X}, i_{Y}] & = & 0\ .               
\end{eqnarray}

\begin{exam}\label{exam2}\emph{ Of course, Lie algebras and their cohomology are basic examples. Another 
extreme example is the tangent Lie algebroid $TM$ of a manifold (when $\rho$ is the identity). Then
one recovers the usual DeRham cohomology, while its representations 
are precisely vector bundles over $M$ endowed with a flat connection. Another important class of examples
are the Lie algebroids associated to Poisson manifold (see section \ref{Poiss}). Foliations $(M, \F)$ can be naturally viewed
as Lie algebroids with $\mathfrak{g}= \F$ the vectors tangent to the leaves, and $\rho$
is the inclusion. Its representations are precisely the foliated vector bundles over $(M, \F)$ \cite{KT, Mol},
while its cohomology with constant coefficients  is well known under the name of foliated 
or leafwise cohomology, denoted $H^*(\F)$ (see e.g. \cite{KaTo, Mol, MoSo, Vais}).
The normal bundle $\nu$ of $\F$ a basic example of representation of $\F$ (see below),
and the cohomology $H^*(\F; \nu)$ is known to be relevant to deformations of foliations \cite{Hei}.
The $\alpha$-simply connected integration of $\F$ is usually called the monodromy groupoid of $\F$
(see \cite{CM} for a description of all groupoids integrating $\F$).
Note that when $\F= \F(\pi)$ is the foliation induced by a submersion $\pi: M\rmap B$ with connected fibers,
then the pull-back $\pi^*E$ of any vector bundle over $B$ is naturally a representation of $\F(\pi)$. Moreover,
if the fibers of $\pi$ are simply connected, then any representation of $\pi$ is of this type. Indeed, in this
case the $\alpha$-simply connected integration of $\F(\pi)$ is just $M\times_{B}M$, which is Morita equivalent to $B$ 
(see Examples \ref{exam1}). }
\end{exam}

\begin{exam}\label{Bottex} {\bf (Bott representations)}\emph{ If $E\subset \mathfrak{g}$ is an ideal (i.e. $[\Gamma E, \Gamma\mathfrak{g}]\subset \Gamma E$), on 
which the anchor vanishes, then $L_{X}(V)= [X, V]$ defines an action of $\mathfrak{g}$ on $E$. If $E$ is abelian (i.e. $[\Gamma E, \Gamma E]=0$),
then it factors through an action of $\mathfrak{h}= \mathfrak{g}/E$ on $E$. This applies in particular to the kernel of the anchor map of a regular Poisson manifold (see the last section).\\
\hspace*{.3in}A similar construction applies to the quotient vector bundle $\nu= \mathfrak{g}/\mathfrak{h}$, where $\mathfrak{h}\subset \mathfrak{g}$
is any Lie sub-algebroid (i.e. closed under the bracket). In this case, the Bott-type formula \cite{Bo}
$L_{X}(\bar{Y})= \overline{[X, Y]}$ makes $\nu$ into a representation of $\mathfrak{h}$.}
\end{exam}

\begin{exam}\label{pull-back}\emph{ Let $\pi: P\rmap M$ be a submersion with connected fibers, and let $\mathfrak{g}$ be a Lie algebroid over $M$.
Recall \cite{McHi} that one has an induced pull-back algebroid $\pi^{!}\mathfrak{g}$ over $P$. 
Its  fiber at $p\in P$ consists of pairs $(X, V)$ with $X\in \mathfrak{g}_{\pi(p)}$, $V\in T_pP$
satisfying  $\rho(X)= (d\pi)_{p}(V)$, its anchor is $(X, V)\mapsto V$. To describe the bracket,
we represent the sections of $\Gamma\pi^{!}\mathfrak{g}$ as sums of elements of type $\phi \pi^{*}(X)$,
with $\phi\in C^{\infty}(P)$, $X\in \Gamma\mathfrak{g}$, and we put:
\[ [(\phi \pi^{*}(X), V), (\psi \pi^{*}(Y), W)]= (\phi\psi\pi^{*}([X, Y])+  L_V(\psi)\pi^{*}(Y)-  L_W(\phi)\pi^{*}(X), [V, W]) .\]
(of course one has to use sums here).We leave to the reader to show that, for any $E\in Rep(\mathfrak{g})$, 
the pull-back bundle $\pi^{*}(E)$ is naturally a representation
of $\pi^{!}\mathfrak{g}$. The map $\Gamma\mathfrak{g}\rmap \Gamma\pi^{!}\mathfrak{g}$ induces a homomorphism
$H^{*}(\mathfrak{g}; E)\rmap H^{*}(\pi^{!}\mathfrak{g}; \pi^{*}E)$ for any representation $E$. 
Note that if $\pi$ has simply connected fibers, then we have an isomorphism $\pi^{*}: Rep(\mathfrak{g})\cong Rep(\pi^{!}\mathfrak{g})$ (use the case mentioned at the end of Examples \ref{exam2}, when $\mathfrak{g}= 0$ and $\pi^{!}(0)= \F(\pi)$ 
is the foliation defined by the fibers of $\pi$). This is only one of the reasons for which any notion of Morita equivalence for Lie algebroids
should declare $\mathfrak{g}$ and $\pi^{!}\mathfrak{g}$ to be equivalent. Hence the following lemma, which will be useful later,
also shows that Morita invariance of cohomology can be expected under certain $n$-connectedness condition.}
\end{exam}

\begin{st}\label{fibers} If $\pi: P\rmap M$ is a submersion with homologically $n$-connected fibers, then
$\pi^*: H^{k}(\mathfrak{g})\rmap H^{k}(\pi^{!}\mathfrak{g})$ is an isomorphism in all degrees
$k\leq n$. The same holds with general coefficients.
\end{st}

{\it Proof:} We first assume that $\mathfrak{g}= 0$; then $\pi^{!}\mathfrak{g}$ is the foliation $\F(\pi)$ defined by the fibers of $\pi$. 
In general, for any foliation $\F$ on $P$, $H^{*}(\F)$ is isomorphic
to the cohomology of $P$ with coefficients on the sheaf $\A$ of smooth functions on $P$ which are locally
constant along leaves. Indeed, $U\mapsto C^{*}(\F|_{U})$ are fine sheaves,
and the resulting complex of sheaves is a resolution of $\A$, as can be seen by restricting to
foliation charts; for more details, see e.g. \cite{MoSo}. In our case, $\A$ is precisely the pull-back of the
sheaf of smooth functions on $M$, and the result is a very special case of a known criterion in sheaf theory
(see e.g. $1.9.4$ in \cite{BeLu}).\\
\hspace*{.3in}In general, we need the spectral sequence associated to
a sub-algebroid $\mathfrak{h}\subset \mathfrak{g}$. At  the first level, it is
\[   H^{p}(\mathfrak{h}; \Lambda^{q}\nu^{*}) \Longrightarrow H^{p+q}(\mathfrak{g}) \ .\]
Here $\nu= \mathfrak{g}/\mathfrak{h}$ is as in Examples \ref{Bottex}. 
This extends the standard spectral sequence for Lie algebras, and the well-known foliated
spectral sequence \cite{KaTo} (when $\mathfrak{h}= \F$ is a foliation, and $\mathfrak{g}= TM$).
To construct it, we consider the filtration $F_{*}C^{*}$ of $C^{*}(\mathfrak{g})$ with $F_{0}C^{*}= C^{*}(\mathfrak{g})$, and, for $q\geq 1$,
\[ F_{q}C^{n}= \{ \omega\in C^{n}(\mathfrak{g}): \omega(X_1, \ldots, X_n)= 0 \ {\rm if} \ X_1, \ldots , X_{n-q+1}\in \Gamma\mathfrak{h} \} .\]
At the $0^{{\rm th}}$ level we clearly have $F_{q}C^{n}/ F_{q+1}C^{n} \cong C^{n-q}(\mathfrak{g}; \Lambda^{q}\nu^{*})$, and, a short computation
shows that the boundary is precisely the one computing $H^{*}(\mathfrak{g}; \Lambda^{q}\nu^{*})$.\\
\hspace*{.3in}With the notations of the statement, $\F(\pi)$ is obviously a Lie sub-algebroid of $\pi^{!}\mathfrak{g}$,
and we consider the associated spectral sequence
\[ E^{p, q}_{1}= H^{p}(\F(\pi); \pi^{*}(\Lambda^{q}\mathfrak{g}^{*})) \Longrightarrow H^{p+q}(\pi^{!}\mathfrak{g})\ .\]
By the case $\mathfrak{g}= 0$ we know that $E^{p, q}= 0$ for $1\leq p\leq n$, hence $H^{k}(\pi^{!}\mathfrak{g})$ in degrees $k\leq n$ 
is isomorphic to the cohomology of the complex $(H^{0}(\F(\pi); \pi^{*}(\Lambda^{*}\mathfrak{g}^{*})), d_{1})$, which is nothing but
$C^{*}(\mathfrak{g})$. The rest is standard. \ $\Boxe$\\

\section{A Van Est isomorphism}
\label{AVEi}

\hspace*{.3in}The most general form of the Van-Est {\it isomorphism} states that, if a
connected Lie group $G$ acts properly (from the right) on a contractible manifold $X$, then
\begin{equation}\label{genVE} 
H^{*}_{d}(G) \cong H^{*}_{G-inv}(X) \ ,
\end{equation}
where the right hand side is the cohomology defined by the complex $\Omega^{*}(X)^{G}$
of $G$-invariant forms on $X$ (cf. $5.6$ in \cite{BW}). In the case where $X= G/K$ with $K\subset G$ maximal compact 
subgroup, $\Omega^*(G/K)^{G}= C^{*}(\mathfrak{g}, K)$ is (by definition) the complex computing the relative Lie algebra cohomology, and (\ref{genVE}) takes the well-known form \cite{VEst, vEst, EstCartan, difcoh}
$H^{*}_{d}(G)\cong H^{*}(\mathfrak{g}, K)$. \\

\hspace*{.1in}In this section we improve (\ref{genVE}), extend it to groupoids,
and present some applications.

\subsection{The proper case}

\hspace*{.3in}In this subsection we prove a particular case
of our general Van Est isomorphism, which is an analogue of
the vanishing of $H^{*}_{d}(G)$ for compact Lie groups $G$.
Recall that a groupoid $\G$ is called proper if the
map $(\alpha, \beta): \nG{1}\rmap \nG{0}\times\nG{0}$ is proper.

\begin{prop}\label{properth} For any proper Lie groupoid $\G$, and any $E\in Rep(\G)$,
\[ H_{d}^{k}(\G; E)= 0, \ \ \ \forall\  k\geq 1 .\]
\end{prop}

{\it Proof:} We fix a smooth Haar system $\lambda$ for $\G$ \cite{Re}, i.e. a family
$\lambda= \{ \lambda^x: x\in \nG{0}\}$ of smooth measures $\lambda^x$ 
supported on the manifolds $\G(-, x)= \beta^{-1}(x)$, with the property that:\\
\hspace*{.3in} $(i)$ for any $\phi \in C_{c}^{\infty}(\nG{1})$, the formula
\[ I_{\lambda}(\phi)(x)= \int_{\G(-, x)} \phi(g) d\lambda^x(g) \]
defines a smooth function $I_{\lambda}(\phi)$ on $\nG{0}$;\\
\hspace*{.3in}$(ii)$  $\lambda$ is left invariant, i.e., for any
$g\in \G(x, y)$, and any $\phi\in C_{c}^{\infty}(G(-, y))$, one has
\[ \int_{\G(-, x)} \phi(gh) d\lambda^x(h)= \int_{\G(-, y)} \phi(h) d\lambda^y(h) . \]
The existence of Haar system for Lie groupoids is well known (see e.g. the preliminaries of \cite{CM}, or \cite{WeNi}). The properness of $\G$ 
ensures the existence of a ``cut-off'' function for $\G$, i.e. a smooth function
on $\nG{0}$ satisfying:\\
\hspace*{.3in} $(iii)$ $\beta: supp(c\compose \alpha)\rmap \nG{0}$ is proper;\\
\hspace*{.3in} $(iv)$ $\int_{G(-, x)} c(\alpha(g)) d\lambda^x(g)= 1$ for all $x\in \nG{0}$.\\
(see e.g. the appendix in \cite{Tu} for the construction of such functions). 
We now check that the following formula defines a contraction of $C^{*}_{d}(\G; E)$:
\[ h(\phi)(g_1, \ldots, g_n) = \int_{\G(-, \beta(g_1))} a\cdot \phi(a^{-1}, g_1, \ldots, g_n) c(\alpha(a)) d\lambda^{\beta(g_1)}(a)\ .\]
The integration is defined because of $(iii)$ above, and defines
a smooth section $h(\phi)$, by $(ii)$ above. One has:
\begin{eqnarray}
 h(d(\phi))(g_1, \ldots, g_n) & = & \int_{\G(-, x)} a\cdot \{ a^{-1}\phi(g_1, \ldots, g_n)- \phi(a^{-1}g_1, g_2, \ldots, g_n)-\nonumber\\
                              & - & \sum_{i=1}^{n-1}(-1)^i \phi(a^{-1}, g_1, \ldots, g_{i}g_{i+1}, \ldots, g_n) +\nonumber\\
                              & + & (-1)^{n+1}\phi(a^{-1}, g_1, \ldots g_{n-1}) \} c(\alpha(a)) d\lambda^x(a) \ ,\nonumber
\end{eqnarray}
which, by $(iv)$ above, and the left invariance $(ii)$ equals to
\begin{eqnarray}
 & \phi(g_1,\ldots , g_n)  -  \int_{\G(-, \beta(g_2))} g_1a\phi(a^{-1}, g_2, \ldots, g_n) c(\alpha(a)) d\lambda^{\beta(g_2)}(a) -\nonumber \\
 & - \sum_{i=1}^{n-1} (-1)^i\int_{\G(-, \beta(g_1))} a\phi(a^{-1}, g_1, \ldots, g_{i}g_{i+1}, \ldots, g_n) c(\alpha(a))d\lambda^{\beta(g_1)}(a) +\nonumber\\
 & +   (-1)^{n+1}\int_{\G(-, \beta(g_1))} a\phi(a^{-1}, g_1, \ldots, g_{n_1}) c(\alpha(a))d\lambda^{\beta(g_1)}(a) =\nonumber\\
 & = \phi(g_1,\ldots , g_n)  -  g_1 h(\phi)(g_2, \ldots, g_n)-\nonumber\\
 & -\sum_{i=1}^{n-1}(-1)^{i}  h(\phi)(g_1, \ldots, g_{i}g_{i+1}, \ldots, g_n)  +  (-1)^{n+1}h(\phi)(g_1, \ldots, g_{n-1}) = \nonumber\\
 & = \phi(g_1,\ldots , g_n)  -  d(h(\phi))(g_1, \ldots, g_n) \ ,\nonumber
\end{eqnarray}
hence the desired formula $hd+ dh= id$.\ $\Boxe$\\

\subsection{A Van-Est theorem}
\label{VanEst}

\hspace*{.3in}We now state and prove the extension of the Van Est isomorphism (\ref{genVE}), mentioned at the begining
of this section.\\
\hspace*{.3in}Let $P$ be a right $\G$-space. We call it proper if the map $P*\nG{1}\rmap P\times P$, $(g, p)\mapsto (gp, p)$ is proper.
Note that, if the moment map $\pi: P\rmap \nG{0}$ of the action is a submersion, then $\G$ acts on the foliation $\F(\pi)$
induced by $\pi$ on $P$. More precisely, for any $p\in P$, and any arrow $g: x\rmap y$
of $\G$ ending at $y= \pi(p)$, the differential at $p$ of the multiplication
$\pi^{-1}(y)\rmap \pi^{-1}(x)$ by $g$ induces a map
\begin{equation}\label{actfol}
g: \F(\pi)_p \rmap \F(\pi)_{pg} \ .
\end{equation}
Hence it makes sense to talk about the complex $C^{*}_{\G}(\F(\pi))$ of $\G$-invariant foliated forms
(a subcomplex of $C^*(\F(\pi))$); denote by $H^{*}_{\G}(\F(\pi))$ the resulting cohomology. Remark that the product structure
on the foliated cohomology induces a product structure on $H^{*}_{\G}(\F(\pi))$.  Also, if $E\in Rep(\G)$,
then $\pi^*E$ is a representation of $\F(\pi))$ (since it is a pull-back via $\pi$), and there are obvious versions 
with coefficients, for $C^{*}_{\G}(\F(\pi); E)$ and $H^{*}_{\G}(\F(\pi); E)$.

\begin{st}\label{VEth} Let $\G$ be a Lie groupoid, and let $P$ be a proper $\G$-space whose moment map $\pi: P\rmap \nG{0}$ is a 
submersion with connected fibers. For any $E\in Rep(\G)$ there is a map compatible with the product structure
\[ \Phi_{P}: H^{*}_{d}(\G; E) \rmap H^{*}_{\G}(\F(\pi); E) .\]
\hspace*{.3in}Moreover, if the fibers of $\pi$
are homologically $n$-connected (i.e. have trivial cohomology in degrees $\leq n$), then
$\Phi_P$ is an isomorphism in degrees $* \leq n$, and is injective
in degree $n+1$.\\
\hspace*{.3in}In particular, if $\pi$ has contractible fibers, then $\Phi_P$ is
an isomorphism in all degrees.
\end{st}

{\it Proof:} We assume that the coefficients are trivial (in general there are obvious modifications).
Consider the space $P*\nG{p}$ consisting of pairs $(p, \overrightarrow{g})$ with $p\in P$, $\overrightarrow{g}\in \nG{p}$ an arrow
of type (\ref{string}) with $\pi(p)= x$. This space comes equipped with a foliation $\F(p)$, defined 
by the projection into $\nG{p}$. In particular, $\F(0)= \F(\pi)$.
We form a double complex $C$ 
which, in bi-degree $(p, q)$, is
\[ C^{p, q}= C^q(\F(p))= C^{\infty}(P*\nG{p}; \Lambda^{q}\F(p)^{*}) .\]
We now describe the differentials (but we leave to the reader the lengthy but straightforward
verification of the compatibility between the horizontal and the vertical differentials). \\ 
\hspace*{.3in}$1$. Columns: the $p^{th}$ column is just the complex $C^*(\F(p))$ computing the foliated cohomology.
Note that it comes with $C^{p}_{d}(\G)$ as a co-augmentation (as the kernel of the differential
$C^0(\F(p))\rmap C^1(\F(p))$):
\begin{equation}\label{rows} 
0\rmap C^{\infty}(\nG{p})\stackrel{\epsilon}{\rmap} C^{0}(\F(p)) \rmap C^{1}(\F(p)) \rmap \ldots 
\end{equation}
\hspace*{.3in}$2$. Rows: the $q^{th}$ row comes with $C^{q}_{\G}(\F(\pi))$ as a co-augmentation:
\begin{equation}\label{columns} 
C^{\infty}(P; \Lambda^q\F(0)^{*})^{\G} \stackrel{\eta}{\rmap} C^{\infty}(P; \Lambda^q\F(0)^{*}) \rmap 
C^{\infty}(P* \nG{1}; \Lambda^q\ \F(1)^{*})  \rmap \ldots   
\end{equation}
To define this, we consider the crossed product groupoid $P\cross \G$. The action
(\ref{actfol}) of $\G$ on $\F(\pi)= \F(0)$ translates into the fact that $\F(0)$ (hence also
$\Lambda^{q}\F(0)^{*}$) is a representation of $P\cross \G$. Note also that
$(P\cross \G)^{\ps p\pd}= P* \nG{p}$, and $\F(p)$, as a vector bundle,
is the pull-back of $\F(0)$. Now, the co-augmentation
$\eta$ is the obvious inclusion, while the rest of (\ref{columns})
(i.e. the $q^{th}$ row $C^{*, q}$) is defined as the complex 
$C^{*}_{d}(P\cross \G; \Lambda^q\F(p)^{*})$  computing the differentiable
cohomology of $P\cross \G$ with coefficients.\\

\hspace*{.1in} Since the properness of $P$ as a $\G$-space
is equivalent to $P\cross \G$ being proper, Proposition \ref{properth} implies that
(\ref{columns}) is exact. Hence the inclusion $\eta: C^{*}_{\G}(\F(\pi))\rmap C$ 
induces isomorphisms:
\begin{equation}\label{etaiz} 
\eta: H^{*}_{\G}(\F(\pi))\cong H^*(C)\ .
\end{equation}
Combined with the inclusion $\epsilon: C^{*}_{d}(\G)\rmap C$, this induces the desired map:
\begin{equation}\label{phiP} 
\Phi_{P}: H^{*}_{d}(\G) \stackrel{\epsilon}{\rmap} H^{*}(C) \stackrel{\eta^{-1}}{\rmap} H^{*}_{\G}(\F(\pi))\ .
\end{equation}
\hspace*{.3in}To prove that $\Phi_{P}$ is compatible with the products, note that $C^{p, q}\otimes C^{p\,', q\,'}$ $\rmap$ $C^{p+p\,', q+ q\,'}$,
\begin{eqnarray}
 & (\omega\cdot\eta)(x, g_1, \ldots , g_{p+p\,'}) = \nonumber\\
 & \omega(x, g_1, \ldots , g_p) \wedge (g_1\ldots g_p)^{-1} \eta( xg_1\ldots g_p, g_{p+1}, \ldots , g_{p+ p\,'}) , \nonumber
\end{eqnarray}
defines a product on $C$, and $\epsilon$ and $\eta$ are maps of DG algebras.\\
\hspace*{.3in}For the second part of the theorem, the spectral sequence
of the double complex $C$, combined with (\ref{etaiz}), provides us with a spectral sequence
\[ E^{p, q}_{2}= H^{p}(H^{q}(\F(*))) \Longrightarrow H^{p+q}_{\G}(\F(\pi)),\]
with $E_{2}^{p, 0}= H^{p}_{d}(\G)$, and with $\Phi_{P}$ as edge maps. Since each $\F(p)$ is a foliation
defined by a submersion with homologically $n$-connected fibers, it follows from Theorem \ref{fibers} (applied to the zero algebroid)
that $E^{p, q}_{2}= 0$ for $1\leq q\leq n$, hence the statement follows by the
well known spectral-sequence arguments. \ $\Boxe$\\

\begin{rmk} {\bf (connection with classifying spaces)}\emph{ This result is undoubtelly related to the classifying space $B\G$ of $\G$ \cite{Se}, or, even better, to the classifying space for proper actions, familiar to the people working on the Baum-Connes conjecture. We discuss here the
classifying space only. In general it is only defined up to homotopy, and this is important for choosing explicit models (depending on the context). Relevant for us is that it is the base space of a principal $\G$-bundle $\sigma: E\G\rmap B\G$, whose  moment map $\pi: E\G\rmap \nG{0}$ has  contractible fibers. We say that $\G$ has a smooth classifying space $B\G$ if these choices can be made in the smooth category. 
Note that in this case, there is a natural (classifying?) Lie algebroid $\tilde{\mathfrak{g}}$ over $B\G$, namely the Lie algebroid of the gauge groupoid (see Examples \ref{exam1}) of $E\G$. Alternatively, $\F(\sigma)= \pi^*\mathfrak{g}$ naturally acts on $\F(\pi)$, and $\mathfrak{g}$ is the quotient algebroid. We deduce the following
result, which extends Haefliger's conjecture mentioned in the last part of our Examples \ref{example1}. }
\end{rmk}

\begin{cor}\label{Haeflconj} Let $\G$ be a Lie groupoid with smooth classifying space $B\G$, and let $\tilde{\mathfrak{g}}$ be
the induced Lie algebroid over $B\G$. Then there are isomorphisms compatible with the product structures:
\[   H^{*}_{d}(\G) \cong H^{*}(\tilde{\mathfrak{g}})  .\]
\end{cor}

\subsection{Application: relation between differentiable and algebroid cohomology}
\hspace*{.3in}Let $\G$ be a Lie groupoid. Roughly speaking,
Lie algebroid cocycles  can be viewed as an infinitesimal version
of groupoid cocycles. This translates into the existence of a map,
\[ \Phi: H^{p}_{d}(\G) \rmap H^{p}(\mathfrak{g}), \]
which we call  {\it the Van Est map}. In the case of Lie groups
is was constructed by Van Est \cite{VEst, vEst}, and it was extended to Lie groupoids
by Weinstein and Xu \cite{WeXu}. Explicitly, $\Phi$ is defined at the chain level
by:
\begin{equation}\label{defphi} 
\Phi(c)(X_1, \ldots , X_p)= \sum_{\sigma\in S_p} sign(\sigma) R_{X_{\sigma\ps p\pd}} \ldots R_{X_{\sigma\ps 1\pd}} (c) , 
\end{equation}
for all $c\in C^{p}_{d}(\G)$, $X_1, \ldots, X_p\in \Gamma(\mathfrak{g})$. Here, for $X\in \Gamma(\mathfrak{g})$
and $c\in C^{p}_{d}(\G)$, $R_X(c)\in C^{p-1}_{d}(\G)$ is given by
\begin{equation}\label{derg} 
 R_{X}(c)(g_2, \ldots, g_p))= \frac{\partial c(-, g_2, \ldots, g_p)}{\partial X} (\beta(g_2)) \ ,
\end{equation}
the derivative at the identity
of $c(-,  g_2, \ldots, g_p): \G(\beta(g_2), -)\rmap \mathbb{C}$,
along the vector field on $\G((\beta(g_2), -)$ induced by $X$. \\
\hspace*{.3in}This construction, together with the fact that
it is an isomorphism in degree $p=1$ provided $\G$ is $\alpha$-simply connected,
are part of the main results/tools of \cite{WeXu}. In the same paper it is conjectured
that, if $\G$ has $2$-connected $\alpha$-fibers, then $\Phi$ is an isomorphism in 
degree $p=2$. In this section we explain how the Van Est isomorphism of the previous 
section clarifies the connection between differentiable and algebroid cohomology (see the theorem below).\\
\hspace*{.3in}Let us first remark that there is a version of $\Phi$ in the presence of coefficients $E$:
one only has 
to replace $c(-, g_2, \ldots , g_p)$ in formula (\ref{derg})
by the function $G(\beta(g_2), -)\rmap E_{\beta(g_2)}$, $g\mapsto g^{-1}c(g, g_2, \ldots, g_p)$.

\begin{st}\label{theWeXu} Let $\G$ be an $\alpha$-connected Lie groupoid, and let $E$ be
a representation of $\G$. The formula (\ref{defphi}) defines a map
\[ \Phi: H^{p}_{d}(\G; E) \rmap H^{p}(\mathfrak{g}; E) \]
which is compatible with the product structures.\\
\hspace*{.3in}Moreover, if the $\alpha$-fibers of $\G$ 
are homologically $n$-connected, then $\Phi$ is an isomorphism 
in degrees $p\leq n$, and is injective for $p= n+1$. 
\end{st}

{\it Proof:} We apply Theorem \ref{VEth} to the right $\G$- space $P= \nG{1}$,
with $\pi= \alpha$ as moment map, and the obvious action. Note that $P$ is always 
proper as a $\G$-space, and $C^{*}_{\G}(\F(\pi); E)$ is isomorphic to $C^{*}_{d}(\mathfrak{g}; E)$
(since $\F(\pi)\cong \beta^*\mathfrak{g}$). 
Hence the only thing we have to prove is that the map $\Phi_{P}$ of Theorem \ref{VEth} coincides
with the map $\Phi$ (\ref{defphi}) of the statement. We use the notations from the proof of Theorem \ref{VEth}, with $P= \nG{1}$, and, as there,
we also assume that the coefficients are trivial (as explained above, the difference with
the general case is mainly notational). Since
$\F(p)\cong \beta^{*}_{p+1}\mathfrak{g}$ (see (\ref{betamap})), $C$ stands for the double complex
$C^{p, q}= C^{\infty}(\nG{p+1}; \Lambda^q(\beta^*\mathfrak{g}^*))$, and we view its elements
as $C^{\infty}(\nG{0})$ multilinear maps
\[ \Gamma(\mathfrak{g})\times \ldots \times\Gamma(\mathfrak{g}) \ni (X_1, \ldots , X_q) \mapsto c(X_1, \ldots , X_q)\in C^{\infty}(\nG{p+1}) .\]
We need the explicit formulas for the horizontal boundary $d^h$ (along $q$), and the vertical one $d^v$ (along $p$) of $C$. First of all
it is not difficult to see that
\begin{equation}\label{dhor}
d^v(c)(X_1, \ldots, X_q)= d^{\,'}(c(X_1, \ldots, X_q))\ ,
\end{equation}
where
\begin{eqnarray} 
 & d^{\,'}: C^{\infty}(\nG{p})\rmap C^{\infty}(\nG{p+1}), \nonumber \\
 & d^{\,'}(c)(g_1, \ldots g_{p+1}) = \sum_{i=1}^{p} (-1)^i c(g_{1}, \ldots , g_{i}g_{i+1}, \ldots , g_{p+1})+ (-1)^{p+1}c(g_{1}, . . . , g_{p}) .\nonumber
\end{eqnarray}
Secondly, since $\F(p)\cong \beta^{*}_{p+1}\mathfrak{g}$, we see that any section $X\in \Gamma(\mathfrak{g})$ induces a vector field
$\tilde{X}$ on $\nG{p+1}$ which is tangent to $\F(p)$. Explicitly, the flow of $\tilde{X}$ is given by:
\begin{equation}\label{flowtilde}
\phi_{\tilde{X}}(g_1, \ldots , g_{p+1})= (\phi_X(t, g_1), g_2, \ldots, g_{p+1}) .
\end{equation}
Since $d^h$ is the boundary
of the complex $C(\F(p))$, we see that the formula for $d^h$ is the usual one: 
\begin{eqnarray}\label{dver}
d^h(c)(X_1, \ldots , X_{q+1}) & = & \sum_{i<j}
(-1)^{i+j-1}\omega([X_i, X_j], X_1, \ldots , \hat{X_i}, \ldots ,
\hat{X_j}, \ldots X_{q+1})) \nonumber \\
 & + & \sum_{i=1}^{q+1}(-1)^{i}
L_{X_i}(\omega(X_1, \ldots, \hat{X_i}, \ldots , X_{q+1})), 
\end{eqnarray}
with the warning that here $L_{X}$ stands for the derivation on $C^{\infty}(\nG{p+1})$ determined by $\tilde{X}$. Note also that, in terms
of these derivations,
\begin{equation}\label{Rder}
 R_{X}(c)= L_{X}(c)\compose u,\ \ \ \ \forall c\in C^{\infty}(\nG{p+1}),
\end{equation}
where
\[ u: \nG{p}\rmap \nG{p+1},\ \ \ u(g_{1}, \ldots, g_{p})= (\beta(g_1), g_{1}, \ldots, g_{p}) .\]
Looking at the construction (\ref{phiP}) of $\Phi_{P}$, to prove that it coincides (in cohomology) with $\Phi$,
it suffices to construct a chain map $\tilde{\Phi}: C\rmap C(\mathfrak{g})$ such that
$\tilde{\Phi}\compose \epsilon= \Phi_{P}$, and $\eta\compose \tilde{\Phi}= id$. 
For
$c\in C^{p, q}$, we define $\tilde{\Phi}(c)\in C^{p+q}(\mathfrak{g})$ by
\begin{eqnarray} 
& \tilde{\Phi}(c)(X_1, \ldots, X_{p+q}) = \nonumber \\
& \sum_{\sigma\in S(p, q)} sign(\sigma) \Phi( c(X_{\sigma\ps 1\pd}, \ldots , X_{\sigma\ps q\pd})\compose u)
(X_{\sigma\ps q+1\pd}, \ldots, X_{\sigma\ps q+ p\pd}) ,\nonumber
\end{eqnarray}
where the sum is over all $(p, q)$- shuffles. E.g., for $c\in C^{p, 1}$,
\begin{equation}\label{particular} 
\tilde{\Phi}(c)(X_1, \ldots, X_{p+1})= \sum_{i} (-1)^{i+1} \Phi( c(X_i)\compose u)(X_1, \ldots, \hat{X_i}, \ldots , X_{p+1}) ,
\end{equation}
and, for $c\in C^{p, 2}$, 
\begin{equation}\label{particulardoi}
\tilde{\Phi}(c)(X_1, \ldots, X_{p+2})= \sum_{i<j} (-1)^{i+j+1} \Phi(c(X_i, X_j)\compose u)(X_1, \ldots, \hat{X_i}, \ldots , \hat{X_j}, \ldots, X_{p+2}).
\end{equation}
The only serious challenge is to prove that
$\tilde{\Phi}$ defines a chain map from (the total complex of) $C$ into $C(\mathfrak{g})$, i.e. we have to prove that
\begin{equation}\label{toprove}
d(\tilde{\Phi}(c))= (-1)^{q} \tilde{\Phi}(d^v(c))+ \tilde{\Phi}(d^h(c)),\ \ \forall\ c\in C^{p, q} .
\end{equation}
The proof of (\ref{toprove}) is actually a combinatorial computation which is based only 
on the following properties of the map $\Phi$:\\
\hspace*{.3in} $(p1).$ $\Phi: C_{d}(\G)\rmap C(\mathfrak{g})$ is a chain map;\\
\hspace*{.3in} $(p2).$ $\Phi(c)= 0$ if $c(g_1, \ldots , g_p)$ does not depend on the variable $g_1$. Note that, together with $(p1)$ this implies that 
$d \Phi(c\compose u)= \Phi(d^{\,'}(c)\compose u) - \Phi(c)$ for all $c\in C_{d}(\G)$; \\
\hspace*{.3in} $(p3).$ for all $c\in C^{\infty}(\nG{p+1})$ one has:
\[ \sum_{i} (-1)^{i+1} \Phi( R_{X_i}(c)) (X_1, \ldots, \hat{X_i} , \ldots , X_{p+1})= \Phi(c)(X_1, \ldots, X_{p+1}) .\]
All these properties are easy to check (and have been remarked also in \cite{WeXu}). 
Since the proof of (\ref{toprove}) for $q=1$ contains all the intricacies of the 
general case, to simplify the exposition we restrict ourselves to this particular case. Hence, let $c\in C^{p, 1}$. Then $d(\tilde{\Phi}(c))(X_1, \ldots, X_{p+2})$ equals to
\[ \sum_{i<j} (-1)^{i+j+1}\tilde{\Phi}(c)([X_i, X_j], X_1, \ldots , \hat{X_i}, \ldots ,
\hat{X_j}, \ldots X_{p+2}) + \sum_{i} L_{X_i}(\tilde{\Phi}(c)(X_1, \ldots, \hat{X_i}, \ldots , X_{p+2}) ,\]
and, by the definition of $\tilde{\Phi}$ (see the particular formula (\ref{particular}) above), this is
\begin{eqnarray}
 & \sum_{i<j} (-1)^{i+j+1} \Phi(c([X_i, X_j])\compose u)(X_1, \ldots, \hat{X_i}, \ldots , \hat{X_j}, \ldots, X_{p+2}) + \nonumber \\
& + \sum_{r<i<j}(-1)^{i+j+1}(-1)^{r} \Phi(c(X_r)\compose u)([X_i, X_j], X_1, \ldots , \hat{X_r}, \ldots ,\hat{X_i}, \ldots ,
\hat{X_j}, \ldots X_{p+2}) + \nonumber\\
& + \sum_{i<r<j}(-1)^{i+j+1}(-1)^{r+1} \Phi(c(X_r)\compose u)([X_i, X_j], X_1, \ldots , \hat{X_i}, \ldots ,\hat{X_r}, \ldots ,
\hat{X_j}, \ldots X_{p+2}) + \nonumber\\
& + \sum_{i<j<r}(-1)^{i+j+1}(-1)^{r} \Phi(c(X_r)\compose u)([X_i, X_j], X_1, \ldots , \hat{X_i}, \ldots ,\hat{X_j}, \ldots ,
\hat{X_r}, \ldots X_{p+2}) + \nonumber\\
& + \sum_{r<i} (-1)^{i}(-1)^{r+1} L_{X_i}( \Phi(c(X_r)\compose u)(X_1, \ldots , \hat{X_r}, \ldots ,\hat{X_i}, \ldots , X_{p+2})) + \nonumber\\
& + \sum_{i<r} (-1)^{i}(-1)^{r} L_{X_i}( \Phi(c(X_r)\compose u)(X_1, \ldots , \hat{X_i}, \ldots ,\hat{X_r}, \ldots , X_{p+2}))\nonumber
\end{eqnarray}
Now, fixing $r$, by the definition of the differential of $C(\mathfrak{g})$, the sum (over $i$) of the last five terms equals $(-1)^{r}d(\Phi(c(X_r)\compose u))(X_1, \ldots, \hat{X_i} , \ldots , X_{p+2})$, which, by $(p2)$ and (\ref{dhor}), is $(-1)^r\Phi(d^v(c)(X_r)\compose u)$ $+$ $(-1)^{r+1}\Phi(c(X_r))$ applied to $(X_1, \ldots, \hat{X_r} , \ldots , X_{p+2})$. Hence our big sum equals to
\begin{eqnarray}
 & \sum_{i<j} (-1)^{i+j+1} \Phi(c([X_i, X_j])\compose u)(X_1, \ldots, \hat{X_i}, \ldots , \hat{X_j}, \ldots,  X_{p+2}) + \label{rah1} \\
 & + \sum_{r} (-1)^{r+1} \Phi(c(X_r))(X_1, \ldots, \hat{X_r}, \ldots , \ldots,  X_{p+2}) - \label{rah2} \\
 & + \sum_{r}  (-1)^r \Phi(d^v(c)(X_r)\compose u)(X_1, \ldots, \hat{X_r} , \ldots , X_{p+2}) \nonumber .
\end{eqnarray}
Since the last term is $-\tilde{\Phi}(d^v(c))(X_1, \ldots, X_{p+2})$ (cf. $(p3)$ above), we
are left with showing the connection of the terms (\ref{rah1}) and (\ref{rah2}) with the horizontal boundary of $C$.
Let us compute $\tilde{\Phi}(d^h(c))(X_1, \ldots, X_{p+2})$. By definition (see (\ref{particulardoi}) above), it is
\[ \sum_{i<j} (-1)^{i+j+1} \Phi(d^h(c)(X_i, X_j)\compose u)(X_1, \ldots, \hat{X_i}, \ldots , \hat{X_j}, \ldots, X_{p+2}) .\]
By the formulas (\ref{dver}) and (\ref{Rder}), this is equal to
\begin{eqnarray}
 & = \sum_{i<j} (-1)^{i+j+1} \Phi(c([X_i, X_j])\compose u)(X_1, \ldots, \hat{X_i}, \ldots , \hat{X_j}, \ldots, X_{p+2}) - \nonumber\\
 & - \sum_{i<j} (-1)^{i+j+1} \Phi(R_{X_i}(c(X_j)))(X_1, \ldots, \hat{X_i}, \ldots , \hat{X_j}, \ldots, X_{p+2}) + \label{cac1}\\
 & + \sum_{i<j} (-1)^{i+j+1} \Phi(R_{X_j}(c(X_i)))(X_1, \ldots, \hat{X_i}, \ldots , \hat{X_j}, \ldots, X_{p+2})  \label{cac2} 
\end{eqnarray}
Replacing the pair $(i, j)$ by $(j, i)$ in (\ref{cac2}), we see that the sum of (\ref{cac1}) and (\ref{cac2}) is the sum over $j$ of
\begin{eqnarray} 
 & (-1)^j \{ \sum_{i<j} (-1)^i \Phi(R_{X_i}(c(X_j))(X_1, \ldots, \hat{X_i}, \ldots , \hat{X_j}, \ldots, X_{p+2}) + \nonumber \\
 & + \sum_{i>j} (-1)^{i+1}  \Phi(R_{X_i}(c(X_j))(X_1, \ldots, \hat{X_j}, \ldots , \hat{X_i}, \ldots, X_{p+2})\}. \nonumber
\end{eqnarray}
By $(p3)$ this is equal to
\[ \sum_{j} (-1)^{j+1} \Phi(c(X_j))(X_1, \ldots , \hat{X_j}, \ldots, X_{p+2}), \]
hence $\tilde{\Phi}(d^h(c))(X_1, \ldots, X_{p+2})$ is precisely the sum of the two expressions (\ref{rah1}) and (\ref{rah2}).
In conclusion, we get the desired relation $d(\tilde{\Phi}(c))= \tilde{\Phi}(d^h(c))- \tilde{\Phi}(d^v(c))$.\ $\Boxe$\\

\begin{rmk}\label{rkfor}\emph{ The formulas (\ref{flowtilde}) and (\ref{Rder}) in the
previous proof show that the map $\Phi$ can be expressed in terms of the flows (see subsection \ref{FGA})
of sections of $\mathfrak{g}$. More precisely, for all $X_1, \ldots, X_p$ $\in\Gamma(\mathfrak{g})$, and $c\in C^{p}_{d}(\G)$,
$R_{X_p} \ldots R_{X_1}(c)(x)$ equals to the derivative at $t_1= \ldots= t_p= 0$ of}
\begin{eqnarray} 
c\mathbb{(}\phi_{X_1}(t_1, \beta_{X_2}(t_2, \beta_{X_3}(t_3, \ldots, \beta_{X_p}(t_p, x) \ldots   )))\ , & & \nonumber \\
                \phi_{X_2}(t_2, \beta_{X_3}(t_3, \ldots, \beta_{X_p}(t_p, x)  \ldots ))\ \ , & & \nonumber \\
                \ \ \ \ \ \ \ \ .  \ \ \ \ \ \ \ \  .         \ \ \ \ \ \ \ \  .     \ \  \ \ \ \ \ \ ,  & & \nonumber \\
                                                         \phi_{X_p}(t_p, x))  \ \ \ \ \ \ \ \ \nonumber
\end{eqnarray}
\emph{In the presence of coefficients $E$, one has to transport the previous element into the fiber $E_x$, i.e. to multiply it by the element}
\[ \phi_{X_1}(t_1, \phi_{X_2}(t_2, \phi_{X_3}(t_3, \ldots, \phi_{X_p}(t_p, x) \ldots\ )))^{-1} .\]
\end{rmk}
 
\begin{rmk}\emph{ Note that, as in the case of classical ``van-Est isomorphisms'' \cite{VEst, vEst, EstCartan}, it is rather the method 
that gives us the maximal information. For instance, the previous double complex induces a spectral sequence which converges
to the algebroid cohomology, and the theorem above is just a simple consequence. By the usual arguments of algebraic topology,
if the fibers of $\G$ are (homological) $n+1$-spheres we obtain a Gysin-type long exact sequence relating differentiable and algebroid cohomology.
Note also that the lower non-trivial part of the spectral sequence provides us with a description of the image of the Van Est map in degree $n+1$. 
To explain this, we consider $n+1$ $\alpha$-loops in $\G$ i.e. smooth
maps $\gamma: S^{n+1}\rmap \G$ with the property that $\alpha\compose \gamma= x$ is constant. For any such $\gamma$ we have a map $\int_{\gamma}: H^{n+1}(\mathfrak{g})\rmap \mathbb{R}$, $[\omega] \mapsto \int_{S^{n+1}} \gamma^{*}\omega|_{L}$, where $L= \G(x, -)$. Note that $\int_{\gamma}$
depends only on the homotopy class of $\gamma$ inside $L$. In particular, if the homotopy class of $\gamma$ is an element of finite order,
then $\int_{\gamma}= 0$. The conclusion is:}
\end{rmk}

\begin{cor} If the $\alpha$-fibers of $\G$ are (homologically) $n$-connected, then $\omega\in H^{n+1}(\mathfrak{g})$ is in the image of 
the Van Est map if and only if $\int_{\gamma}\omega= 0$ for all $n+1$ $\alpha$-loops $\gamma$.
\end{cor}

\begin{rmk}{\bf (compact supports)}\emph{ In relation with the cyclic cohomology of the convolution algebra, it is relevant 
to have a compactly supported version of differentiable cohomology and of our results. We prefer to call this cohomology with
compact supports a homology theory, and denote it by $H^{d}_{*}(\G)$ instead of $H^{*}_{cpt, d}(\G)$; this is in agreement also with
the homology theory described in \cite{CrMo}, which we recover when $\G$ is \'etale. The definition of $H^{d}_{*}(\G)$, in analogy with
the definition of the convolution product \cite{Re}, depends on the choice of a smooth Haar system $\lambda$ for $\G$. For the definition, 
see the proof of our Proposition \ref{properth}. The (chain) complex $C_{*}^{d}(\G)$ defining this homology consists of
the differentiable cochains with compact supports, while the differential is given by}
\begin{eqnarray} 
(dc)(g_{1}, \ldots , g_{p}) & = & \int c(a^{-1}, g_{1}, \ldots , g_{p})d\lambda^{\beta({g_1})}(a) + \\
   & + & \sum_{i=1}^{p} (-1)^i \int c(g_1, \ldots , g_ia, a^{-1}, g_{i+1}, \ldots , g_{p}) d\lambda^{\alpha(g_{i})}(a) + \\
   & + & (-1)^{p+1} \int c(g_1, \ldots , g_{p}, a) d\lambda^{\alpha(g_{p})}(a) .
\end{eqnarray}
\emph{A compactly version of our proofs shows that:\\
\hspace*{.3in}$1.$  $H^{d}_{*}(\G)$ is Morita invariant. In particular (compare to the similar statement for convolution algebras \cite{MRW}), it does not depend
on the choice of the smooth Haar system $\lambda$.\\
\hspace*{.3in}$2.$ Denote by $q$ the dimension of $\mathfrak{g}$, and by $H^{*}_{cpt}(\mathfrak{g})$ the algebroid cohomology with compact supports. 
Assume we are in the orientable case (in general, some twisting in the algebroid cohomology is necessary). 
If $\G$ has homologically $n$-connected fibers, then one has isomorphisms $H^{d}_{k}(\G) \cong H^{q- k}_{cpt}(\G)$ for all $0\leq k\leq n$.}
\end{rmk}

\subsection{Application: an integrability result}
\label{Air}

\hspace*{.3in}As an application of Theorem \ref{theWeXu}, we show the following integrability result. In the case 
of Lie algebras this result is precisely the argument given by Van Est \cite{VEst} for a short proof of Lie's third theorem. 
An immediate consequence will be a more conceptual proof and a slight improvement of 
Dazord-Hector's integrability criterion for Poisson manifolds (this will be explained in the last section).

\begin{st}\label{integration} Let 
\begin{equation}\label{extens} 
0\rmap \mathfrak{l}\stackrel{i}{\rmap} \mathfrak{h} \stackrel{\pi}{\rmap}\mathfrak{g} \rmap 0 
\end{equation}
be an exact sequence of Lie algebroids, with $\mathfrak{l}$ abelian.
If $\mathfrak{g}$ admits a Hausdorff integration whose $\alpha$-fibers 
are simply connected and have vanishing second cohomology groups,
then $\mathfrak{h}$ is integrable 
(and the integration can be chosen to be Hausdorff and $\alpha$-simply connected).
\end{st}

{\it Proof:} We need the following remarks:\\
\hspace*{.3in}$1.$ given a representation $E$ of a Lie algebroid $\mathfrak{g}$ over $M$,
and $\omega\in C^2(\mathfrak{g}; E)$ closed, one can form a twisted semi-direct product Lie algebroid 
$\mathfrak{g}\ltimes_{\,\omega}E$. As a vector bundle it is $\mathfrak{g}\oplus E$, with the anchor
$(X, V)\mapsto \rho(X)$, and with the bracket
\[ [(X, V), (Y, W)]= ( [X, Y], L_{X}(W)- L_{Y}(V)+ \omega(X, Y)] \ ,\]
for $X, Y\in \Gamma(\mathfrak{g}), V, W\in \Gamma(E)$. Moreover, the isomorphism class of 
$\mathfrak{g}\ltimes_{\,\omega}E$ depends only on the cohomology class $[\omega]\in H^{2}(\mathfrak{g}; E)$.
In the trivial case ($\omega= 0$), one obtains the usual semi-direct product $\mathfrak{g}\ltimes E$.
These remarks are well-known \cite{McK} and follow by direct computation.\\
\hspace*{.3in}$2.$ Similarly, given a representation $E$ of a Lie groupoid $\G$ over $M$, and $c\in C^{2}_{d}(\G; E)$ closed,
one forms the groupoid $\G\ltimes_{\,c}E$ over $M$, whose space of objects is $\nG{1}\times_{\nG{0}}E$,
with the product
\[(g_1, v_1) (g_2, v_2)= (g_1g_2, g_{1}^{-1}\cdot v_1+ v_2+ (g_1g_2)^{-1}\cdot c(g_1, g_2))\ .\]
Moreover, $Lie(\G\ltimes_{\,c}E)\cong \mathfrak{g}\ltimes_{\,\omega}E$, where $\omega= \Phi(c)$. Since this is 
obvious as an isomorphism of vector bundles, the more serious task is to identify the bracket of
$Lie(\G\ltimes_{\,c}E)$. For this, one notes that the flows $\phi_{(X, V)}$ ($X\in \Gamma(\mathfrak{g})$, $V\in \Gamma(E)$) for
$\G\ltimes_{\,c}E$ are given by 
\[   \phi_{X, V}(t, x)= (\phi_{X}(t, x), tV(x)),\ \ \ x\in M\ ,\]
and then one uses the formula (\ref{Liebracket}) for the Lie bracket in $Lie(\G\ltimes_{\,c}E)$,
the formula (\ref{actiong}) for the action of $\mathfrak{g}$ on $E$, and the formula described in Remark \ref{rkfor} for the map $\Phi$.\\
\hspace*{.3in}$3.$ Given an extension (\ref{extens}) as in the statement,
there is a canonical action of $\mathfrak{g}$ on $\mathfrak{l}$, and a canonical
cohomology class $[\omega]\in H^{2}(\mathfrak{g}; \mathfrak{l})$ such that $\mathfrak{h}\cong \mathfrak{g}\ltimes_{\,\omega}\mathfrak{l}$.
They are determined by the formulas 
\[ i(L_X(V))= [\sigma(X), i(V)],\ \ i(\omega(X, Y))= \sigma([X, Y])- [\sigma(X), \sigma(Y)]\ ,\]
for $X, Y\in \Gamma(\mathfrak{g}), V\in \Gamma(\mathfrak{h})$, where $\sigma: \mathfrak{g}\rmap \mathfrak{h}$ is a/any linear splitting of $\pi$.
This is an easy computation, whose conclusion is the well known correspondence \cite{McK} between $2$-cohomology classes $[\omega]\in H^2(\mathfrak{g}; \mathfrak{l})$ and extensions (\ref{extens}).\\

\hspace*{.1in}The proof of the theorem is quite obvious now: let $\G$ be the integration
of $\mathfrak{g}$ whose existence is part of the hypothesis. Then, the canonical action of $\mathfrak{g}$ on $\mathfrak{l}$
(cf. $3.$ above) comes from an action of $\G$ on $\mathfrak{l}$  since $\G$ is $\alpha$-simply connected (see subsection \ref{ss1.3}). By theorem \ref{theWeXu},
the canonical $2$-cohomology class $\omega\in H^2(\mathfrak{g}; \mathfrak{l})$ (cf. $3.$ above) is of type $\omega= \Phi(c)$
for some $c\in H^{2}_{d}(\G; \mathfrak{l})$. Then, using $2.$ above, we see that $\Ha= \G\ltimes_{\,c}\mathfrak{l}$ is a
(Hausdorff) integration of $\mathfrak{h}$.\ $\Boxe$\\

\begin{rmk}\label{commint}\emph{ Actually we have proven a bit more then stated:
 given an extension (\ref{extens}), it defines
a cohomology class in $\omega \in H^{2}(\mathfrak{g}; \mathfrak{l})$, and, if this class is integrable 
(i.e. in the image of the Van Est map), 
then $\mathfrak{h}$ is integrable. As explained in the introduction, this theorem immediately
implies Lie's third theorem (see also section 14 in \cite{VEst}). Essential for this is that simply connected
Lie groups have vanishing second cohomology groups, and that we can use the adjoint representation.
For Lie groupoids these both fail. Note that the most dificult problem to overcome is the non-existence of the adjoint representation 
for general Lie algebroids. However, there are no such difficulties in the case of Lie algebra bundles (the LAB's of \cite{McK}),
or, more generally, in the case of bundles of Lie algebras with regular center (i.e. with the property that the centers of the Lie algebra fibers
fit into a vector bundle). The resulting integrability results are  particular cases of Douady-Lazard's theorem \cite{DoLa}. }
\end{rmk}

\section{Characteristic classes in algebroid cohomology}

\hspace*{.3in}The aim of this section is to find the characteristic
invariants living in algebroid cohomology, as well as their relation to the
Van Est map. As we already mentioned, many examples of Lie algebroids arise naturally and, even if they
are integrable, the integrating groupoids may be difficult to visualise (the best examples
this remark applies to are probably the Lie algebroids associated to Poisson manifolds, but this type of phenomenon
occurs even in the case of foliations). For this reason we insist on working at the algebroid level.\\

\hspace*{.1in}Asking for characteristic invariants living in algebroid cohomology,
there is an obvious (and naive) construction. Using the composition with the anchor,
$\rho_{*}: H^*(M)\rmap H^*(\mathfrak{g})$, we define the $\mathfrak{g}$-Chern classes of a vector bundle $E$ over $M$ by
\[ Ch^{\mathfrak{g}}_{k}(E):\,= \rho_{*}( Ch_{k}(E)) \in H^{2k}(\mathfrak{g}),\ k\geq 0.\]
In the case of foliations, these characteristic classes were studied by Moore Schochet
\cite{MoSo} under the name of the foliated Chern classes. They also appear in the geometry of Poisson manifolds $P$, when $E$ is e.g. a co-foliation on $P$ \cite{Vais}. In subsection \ref{uk's} we will prove 
that $Ch^{\mathfrak{g}}_{k}(E)$ can be computed using $\mathfrak{g}$-connections on $E$ 
\cite{ELW}; this will be done in the proof of the following theorem, which shows that
these classes can also be viewed 
as obstructions to the existence of infinitesimal actions of $\mathfrak{g}$ on $E$.

\begin{st}\label{vanishing} For any representation $E$ of $\mathfrak{g}$,
the Chern classes $Ch_{k}^{\mathfrak{g}}(E)\in H^{2k}(\mathfrak{g})$ 
vanish for $k\geq 1$.
\end{st}

\hspace*{.1in}According to the general principle that a vanishing result
for certain characteristic classes is the origin of new (secondary) characteristic 
classes which are made out of the transgression of the vanishing ones, the previous 
proposition rises the question of finding the non-trivial cohomological invariants 
of representations of $\mathfrak{g}$.\\
\hspace*{.3in}In the next subsection we will show in an 
explicit fashion how such a first class arises naturally; the complete description 
of the secondary classes (and the proof of the theorem above) will be given in the 
subsection \ref{uk's}. In \ref{reltoVE} we will describe the relation with 
the Van Est map. \\

{\bf Warning:} So far, the choice of the basic field (real or complex numbers) was irrelevant.
From now on we use the notations $Rep_{\mathbb{R}}(\mathfrak{g})$ and $Rep_{\mathbb{C}}(\mathfrak{g})$
to distinguish between real/complex representations. Also we consider {\it only cohomology with real 
coefficients} (hence $H^{*}(\mathfrak{g})$ will always stand for $H^{*}(\mathfrak{g}; \mathbb{R})$).

\subsection{The first characteristic class $u_1$}
\label{firstclass}

\hspace*{.3in}We now describe, in a direct manner, the first cohomology class
$u_1(E)\in H^1(\mathfrak{g})$ of a representation $E$ of the Lie algebroid $\mathfrak{g}$.
Denote by $M$ the base manifold. We consider the real/complex cases separately.\\

{\bf The real case:} Let $E\in Rep_{\mathbb{R}}(\mathfrak{g})$. We first assume that $E$, as a vector bundle over $M$, is trivializable,
and let $e= \{e_1, \ldots , e_n\}$ be a frame of $E$. Then the equation:
\[ L_X(e_j)= \sum \omega_{j}^{i}(X) e_i, \ \ X\in \Gamma(\mathfrak{g}) \]
shows that the action of $\mathfrak{g}$ on $E$ is uniquely determined by 
a matrix:
\[ \omega_e = (\omega_{j}^{i})_{1\leq i, j\leq n} \in M_n(C^1(\mathfrak{g})) \]
satisfying the flatness condition
\[ d\omega_e= \omega_e\wedge \omega_e= \frac{1}{2}[\omega_e, \omega_e] \ .\]
The last relation implies that
\[ Tr(\omega_e)\in C^1(\mathfrak{g}) \]
is closed. The same happens with the forms
\begin{equation}\label{highertr} 
Tr(\underbrace{\omega_e\wedge \ldots \wedge \omega_e}_{2k-1})\in C^{2k-1}(\mathfrak{g}) ,
\end{equation}
but at these classes we will look in the next subsection.\\
\hspace*{.3in}If $f= \{f_1, \ldots , f_n\}$ is another frame, and $A= (a^{i}_{j})\in M_n(C^{\infty}(M))$ 
is the change of coordinates matrix, i.e. $e_j= \sum a_{i}^{j} f_i$, one can easily see that 
the matrix corresponding to $f$ is
\begin{equation}\label{changebase} \omega_f= A\omega_eA^{-1} + (dA)A^{-1} .
\end{equation}
This implies
\[ Tr(\omega_f)= Tr(\omega_e) + d(log|detA|) ,\]
hence the class:
\begin{equation}\label{trivial}
u_1(E):= [Tr(\omega_e)]\in H^{1}(\mathfrak{g}) 
\end{equation}
does not depend on the choice of the trivialisation.\\

\hspace*{.1in}We now return to the general case, where $E$ is not necessarily
trivializable. We choose a covering $\U = \{U_{\alpha}\}$ of domains of trivialization 
of $E$, and let
\[ h_{\alpha, \beta}: U_{\alpha}\cap U_{\beta} \rmap GL_n\]
be the transition functions of $E$. The previous construction provides us with 
$1$-forms of the restrictions of $\mathfrak{g}$ to $U_{\alpha}$'s:
\[ u(\alpha)= Tr(\omega_{\alpha}) \in C^{1}(U_{\alpha}; \mathfrak{g}) ,\]
and with functions
\[ u(\alpha, \beta)= log|det (h(\alpha, \beta))| \in C^{0}(U_{\alpha}\cap U_{\beta}; \mathfrak{g}) \]
satisfying:
\begin{eqnarray}
 d(u(\alpha)) & = & 0 \ \ {\rm on} \ U_{\alpha}\ ,\nonumber\\
 d(u(\alpha, \beta)) & =  & u(\alpha)- u(\beta) \ \ {\rm on} \ U_{\alpha}\cap U_{\beta}\ ,\nonumber\\
 u(\alpha, \beta)- u(\alpha, \gamma) + u(\beta, \gamma) & = & 0\ \ {\rm on} \ U_{\alpha}\cap U_{\beta}\cap U_{\gamma}\ ,\nonumber
\end{eqnarray}
The last relation follows from the cocycle relations satisfied by the transition functions. These relations
are precisely the relations defining a closed cocycle $u$ in the Cech double complex $\check{C}^{*}(\U; C^*(\mathfrak{g}))$.
By a Mayer-Vietoris argument (cf. e.g. \cite{BoTu}), one has exact sequences 
\[ 0\rmap C^p(\mathfrak{g})\stackrel{r_{\U}}{\rmap} \check{C}^{0}(\U; C^p(\mathfrak{g}))\stackrel{\delta}{\rmap}
\check{C}^{1}(\U; C^p(\mathfrak{g}) \stackrel{\delta}{\rmap} \ldots \]
where $\delta$ is the Cech boundary, and $r_{\U}$ is the obvious restriction to opens $U\in\U$. Hence the map $r_{\U}$ can be viewed
as a map from $C^*(\mathfrak{g})$ into (the total complex of) $\check{C}^{*}(\U; C^*(\mathfrak{g}))$, which induces an isomorphism
in cohomology, and we define
\[ u_{1}(E)= r_{\U}^{-1}([u])\in H^{1}(\mathfrak{g}) .\]
By the same arguments as in the trivializable case:

\begin{prop} Given a representation $E$ of a Lie algebroid $\mathfrak{g}$, the previously constructed class
$u_{1}(E)\in H^{1}(\mathfrak{g})$ does not depend on the choice of the local trivializations, and,
when $E$ is trivializable as a vector bundle, it is given by
the formula (\ref{trivial}).
\end{prop}

{\bf The complex case:} Assume now that $E\in Rep_{\mathbb{C}}(\mathfrak{g})$. We can simply define
$u_1(E)= \frac{1}{2}u_{1}(E_{\mathbb{R}})$, where $E_{\mathbb{R}}$ is the real representation underlying $E$. Although this is
correct, it is quite instructive to try to imitate the real case.
As above, the choice of a trivialization $e$ (over $\mathbb{C}$) 
for $E$ defines a matrix $\omega_{e}$, and
\[ Tr(\omega_e) \in C^{1}(\mathfrak{g}; \mathbb{C}) \]
will be a closed cocycle. However, simple examples show that it is only $Tr(Re(\omega_{e}))$ that is invariant under the change 
of trivializations. Indeed, the analogue of (\ref{changebase}) gives
\[ Tr(\omega_{f})- Tr(\omega_e)= \frac{\alpha d\alpha+ \beta d\beta}{\alpha^2+ \beta^2} + i \frac{\alpha d\beta- \beta d\alpha}{\alpha^2+ \beta^2}  ,\]
where $\alpha+ i \beta= det(A)$. While $(\alpha d\alpha+ \beta d\beta)/(\alpha^2+ \beta^2)$ is $d(log|det(A)|)$, 
\begin{equation}\label{volform}
\frac{\alpha d\beta - \beta d\alpha}{\alpha^2+ \beta^2}
\end{equation} 
is not always exact. Since $Tr(Re(\omega_f))= Tr(Re(\omega_e)) + d(log|detA|)$ we can proceed exactly as in the real case (with $Tr(\omega_e)$
replaced by $Tr(Re(\omega_e))$) and define the class $u_1(E)\in H^1(\mathfrak{g})$. It coincides with $\frac{1}{2}u_{1}(E_{\mathbb{R}})$. \\
\hspace*{.3in}Note that what the complex case really teaches us is that,  if one wants to fully use the closed forms $Tr(\omega_e)$, there is
an obstruction to the independence on $e$, obstruction which comes from the non-exactness of (\ref{volform}). It is worth to point out
that (\ref{volform}) gives the volume form on $S^1= U(1)$. Although this type of phenomenon has not been seen in the real case above, it will be present
in higher dimensions both for real and complex representations. \\

\begin{exam}\label{example2}\emph{ In the case of $1$-dimensional representations $L$,
then our $u_1(L)$ coincides with the characteristic class $\theta_{L}$ introduced in \cite{ELW}. The modular class 
of $\mathfrak{g}$ is defined as $u_1(Q_{\mathfrak{g}})$, where $Q_{\mathfrak{g}}= \Lambda^{top}\mathfrak{g}\otimes \Lambda^{top}T^*M$
is the canonical line bundle of $\mathfrak{g}$. When applied to the Lie algebroid associated to a
Poisson manifold $P$, one obtains (up to a constant) the modular class of $P$. For more details see \cite{ELW}, and also our last section.\\
\hspace*{.3in}Let us now look at a simple example which shows the non-triviality of $u_1$. Let $X\in \xx(M)$ be a vector
field on the manifold $M$. It induces a Lie algebroid $\mathfrak{g}_{X}$ 
over $M$: as a vector bundle it is just $M\times \mathbb{R}$, $\rho$
is the multiplication by $X$, while the Lie bracket on $\Gamma(\mathfrak{g}_X)= C^{\infty}(M)$
is $[f, g]= fX(g)- X(f)g$. For any zero $x$ of $X$, the evaluation at $x$ defines a map
$ev_x: H^1(\mathfrak{g}_{X})\rmap \mathbb{C}$. The cotangent bundle $T^*M$ is a representation 
of $\mathfrak{g}_X$; a local computation shows that $ev_x(u_1(T^*(M)))= \sum_i \frac{\partial X^i}{\partial x^i}(x)$,
where $(x^1, \ldots, x^n)$ is a/any system of local coordinates around $x$.}
\end{exam}

\subsection{The higher characteristic classes $u_{2k-1}$}
\label{uk's}

\hspace*{.3in}The key of finding all the higher cohomological invariants of representations
is a better understanding of the notion of representation, in terms of frame bundles (see Corollary 
\ref{connections} below). This is
analogous to the correspondence between connections on a vector bundle over a manifold,
and connections on the associated principal $GL_n$ (frame) bundle. To any $E\in Rep(\mathfrak{g})$
we will associate certain classes $u_{2k-1}(E)\in H^{2k-1}(\mathfrak{g})$ which are non-trivial only when
$1\leq k\leq dim(E)$, and which for real representation vanish for even $k$'s. To get a feeling about the final result, I mention here
that the forms defining these classes are just some ``corrections'' of (\ref{highertr}) (see also Examples \ref{example3}). \\

{\bf Explicit approach:} (sketch) Let $E$ be a (complex or real) representation of $\mathfrak{g}$. As in the previous 
subsection, we first assume that $E$ has a trivialization $e$. We have noticed that, together with $Tr(\omega_e)$, we also have
the closed forms (\ref{highertr}). Of course, they appear as the natural candidates for the higher cohomological
invariants of $E$. However, as hinted by the complex case above, the resulting classes depend on the trivialization $e$,
hence they cannot be globalized. Instead, we have to ``correct'' these classes (\ref{highertr}) and consider:
\begin{equation}\label{highertrm}
Tr(\underbrace{\theta_e\wedge \ldots \wedge \theta_e}_{2k-1})\in C^{2k-1}(\mathfrak{g}) , \ \ \theta_e= (\omega_e + \omega_{e}^{*})/2 
\end{equation}
(where $\omega^{*}_{e}= (\overline{\omega}_{e})^{t}$).
The full understanding of this choice, and of the fact that these are all geometric classes one can construct, lies in the computation of the differentiable cohomology of $GL_n$, and the global understanding of representations in terms of frame bundles. This will be done 
below. We leave for the reader to show that the previous classes are indeed closed, and that in the real case they vanish if $i$ is even. \\
\hspace*{.3in}For the sake of explicit formulas, I mention here a global expression of these classes (see also \cite{BiLo}). Hence, let $E$ be any representation (real or complex) of $\mathfrak{g}$. Denote by $\nabla$ the action of $\mathfrak{g}$ on
$E$, and by $\nabla^{*}$ be the dual action of $\mathfrak{g}$ on
$E^{*}$. Let $h$ be a metric on $E$. Using the isomorphism $E\cong E^*$ induced by $h$, we transport $\nabla^{*}$ to a new action of 
$\mathfrak{g}$ on $E$, denoted $\nabla^{h}$. Then $\theta(E, h):= (\nabla- \nabla^{h})/2$ is in $C^{1}(\mathfrak{g}; End(E))$,
and, up to a constant, the classes $u_{2k-1}(E)$ will be given by (\ref{highertrm}) with $\theta_{e}$ replaced by 
$\theta(E, h)$. Alternatively, one can use a transgression  (a la Chern-Weil) construction for $(\nabla, \nabla^h)$. 
The resulting cohomology classes will be independent of the choice of the metric. We now turn to the promised approach, in terms of frame bundles.\\

{\bf Global approach:} To obtain all the characteristic classes at once, we restrict to the complex case, and, for $E$ real,
we will define $u_{2k-1}(E)= u_{2k-1}(E_{\mathbb{C}})$, where $E_{\mathbb{C}}$ is the complexification of $E$. Hence 
let $E$ be an $n$-dimensional complex representation of $\mathfrak{g}$. Denote by
$\pi: P\rmap M$ the frame bundle of $E$. We have already seen that, fixing
a frame $\{ e_1, \ldots , e_n\}$ for $E$, the action of $\mathfrak{g}$ on $E$ 
is uniquely determined by
a matrix $\omega_e \in C^{1}(\mathfrak{g})\otimes gl_n$ satisfying $d\omega_e= \frac{1}{2}[\omega_e, \omega_e]$.
On the other hand, the basic properties (\ref{unu}), (\ref{doi}) of an action show (as in the
case of connections on vector bundles) that the expression $L_{X}(s)(x)$, for $x\in M$, depends
only on $X_{x}$ and on the restriction of $s$ to an integral curve of $\rho(X)\in \xx(M)$ through $x$.
Combining the previous two remarks we see that, 
for any $e(0)= \{ e_{1}(0), \ldots, e_{n}(0)\} \in P$ a frame of $E_x$, any $X\in \mathfrak{g}_x$,
and any tangent vector $V\in T_{e(0)}P$ defined by a curve of frames $e: (-\epsilon, \epsilon) \rmap P$ around $e(0)$
such that 
\begin{equation}\label{defpull}
\rho(X)= (d\pi)_{p}(V),
\end{equation}
one has an (unique) associated matrix $\omega_{j}^{i}= \omega_{j}^{i}(X, V)\in gl_n$ 
such that $L_{X}(e_j)= \sum \omega_{j}^{i} e_{i}$.
Now, (\ref{defpull}) is precisely the relation defining the pullback algebroid $\pi^{!}(\mathfrak{g})$ (see Examples \ref{pull-back}). Hence the previous construction provides us with an element
\begin{equation}\label{oneform}
\omega= \omega_E \in C^{1}(\pi^{!}\mathfrak{g}) \otimes gl_{n} .
\end{equation}
To understand the special features of $\omega$ we use 
the canonical Lie algebra map
\begin{equation}\label{canonical}
gl_{n}\rmap \Gamma(\pi^{!}\mathfrak{g}),\ v\mapsto v^{\sharp}= (0, v_{P}) .
\end{equation}
Here $v_{P}\in \xx(P)$ is the transport of $v$ to $P$; it comes
from the differential (at the identity matrix) 
\begin{equation}\label{inc}
j_p: gl_n \rmap T_pP
\end{equation}
of the multiplication map $GL_n \rmap P$, $g \mapsto p\cdot g$. Note that,
since the inclusion
(\ref{inc}) maps $gl_n$ isomorphically into the
space of vertical vector fields on $P$ (i.e. killed by the differential of $\pi$),
 $v^{\sharp}$
does indeed define an element in $\Gamma(\pi^{!}\mathfrak{g})$. 
Now, the Cartan calculus on $\pi^{!}\mathfrak{g}$ (cf. \ref{ss1.3}), 
together with the Lie
algebra map (\ref{canonical}), endows the DG algebra $C(\pi^{!}\mathfrak{g})$
with Lie derivatives (which are derivations of degree $0$) and interior products
(which are derivations of degree $1$),
\[ L_{v}: C(\pi^{!}\mathfrak{g})\rmap C(\pi^{!}\mathfrak{g}),\ i_{v}: C(\pi^{!}\mathfrak{g})\rmap C^{*-1}(\pi^{!}\mathfrak{g}),\]
linear on $v\in gl_n$, and which satisfy
the Cartan relations (\ref{C1})- (\ref{C4}). The Lie derivatives can be viewed also as the infinitesimal 
version of the canonical action of $GL_n$ on $C(\pi^{!}\mathfrak{g})$. 
In the standard terminology which originates in Cartan's interpretation of connections
(see \cite{KT} and the references therein), $C(\pi^{!}\mathfrak{g})$ is
a $gl_n$- DG algebra. The main properties of the $1$-form (\ref{oneform}) are:\\
\hspace*{.3in} $(i)$. $i_v(\omega)= v$ for all $v\in gl_n$,\\
\hspace*{.3in} $(ii)$. $\omega$ is is $GL_n$-invariant,\\
\hspace*{.3in} $(iii)$. $d\omega= \frac{1}{2}[\omega, \omega]$. \\
(i.e., in the terminology of \cite{KT} again, $\omega$ is a flat connection for the $gl_n$-DG algebra $C(\pi^{!}\mathfrak{g})$).
We can interpret  (\ref{oneform}) as a map $gl_{n}^{*}\rmap C^1(\pi^{!}\mathfrak{g})$, which can be uniquely extended to a map of algebras
\begin{equation}\label{kappaE} 
k_{E}: C^*(gl_n) \rmap C^*(\pi^{!}\mathfrak{g}) \ .
\end{equation}
It is easy to see that the basic properties of $\omega$ translate into:\\
\hspace*{.3in} $(i\,')$. $k_E$ is compatible with the interior products $i_v$, $v\in gl_{n}$,\\
\hspace*{.3in} $(ii\,')$. $k_E$ is compatible with the Lie derivatives $L_v$, $v\in gl_{n}$,\\
\hspace*{.3in} $(iii\,')$. $k_E$ is compatible with the differentials. \\
(this is the standard passing from a flat connection $1$-form $\omega$ to a map of $gl_n$-DG algebras).
Moreover, we can also go backwards, hence our discussion can be summarized into:

%
%
%
%

\begin{cor}\label{connections}
Let $\mathfrak{g}$ be a Lie algebroid, let $E$ be a vector bundle over the base manifold $M$ of $\mathfrak{g}$,
and let $\pi: P\rmap M$ be the associated principal $GL_n$-bundle of $E$. Then there is a one to one
correspondence between:\\
\hspace*{.3in}$(1)$ a pairing $\Gamma(\mathfrak{g})\times \Gamma (E)\rmap \Gamma (E)$ which makes $E$ into a representation of $\mathfrak{g}$;\\
\hspace*{.3in}$(2)$ a $1$-form $\omega\in C^1(\pi^{!}\mathfrak{g})\otimes gl_n$ satisfying $(i)$- $(iii)$ above;\\
\hspace*{.3in}$(3)$ a map $k: C^*(gl_n)\rmap C^*(\pi^{!}\mathfrak{g})$ of DG algebras, satisfying $(i\,')$- $(iii\,')$ above.
\end{cor}

\hspace*{.1in}Now, given our representation $E$, $k_{E}$ induces a map at the level of $U(n)$-basic elements, i.e. elements which are
$U(n)$-invariant, and are killed by all $i_{v}$'s with $v\in u_n$,
\begin{equation}\label{Ulrside} 
k_{E}: C^{*}(gl_n)_{U(n)-basic} \rmap C^*(\pi^{!}\mathfrak{g})_{U(n)-basic}\ .
\end{equation}
The reason for passing to $U(n)$-basic elements is that we can get down from $\pi^{!}\mathfrak{g}$ over $P$ to $\mathfrak{g}$ over $M$. Indeed, the right hand side is isomorphic to $C^*(\pi_{0}^{!}\mathfrak{g})$, where $\pi_0: P_0= P/U(n) \rmap M$
is the obvious projection. Its fibers are contractible, hence, by Theorem  \ref{fibers}, it induces isomorphisms
\[ \pi_{0}^{*}: H^*(\mathfrak{g})\tilde{\rmap} H^*(\pi_{0}^{!}\mathfrak{g}) .\]
On the other hand, the right hand side of (\ref{Ulrside}) is precisely the complex computing the relative Lie algebra cohomology $H^{*}(gl_n, U(n))$.
Recall (see e.g. \cite{KT}) that this cohomology is the exterior algebra 
on $n$ generators $u_1, u_3, \ldots , u_{2n-1}$ of degrees $deg(u_{2k-1})= 2k-1$. Hence, from (\ref{Ulrside}) we get a map in cohomology
\begin{equation}\label{thecharactr}  
k_{E}:  \Lambda^*(u_1, u_3, \ldots , u_{2n-1})\cong H^*(gl_n, U(n))\rmap H^*(\mathfrak{g}) .
\end{equation}

\begin{defin} Define the characteristic classes of the $n$-dimensional representation $E$ of the Lie algebroid $\mathfrak{g}$ as:
\[ u_{2k-1}(E):= k_{E}(u_{2k-1})\in H^{2k-1}(\mathfrak{g}), \ \ \ 1\leq i\leq n \ .\]
\end{defin}

\hspace*{-.2in}{\bf Main properties:} The main properties that $u_{2k-1}$'s satisfy are:\\
\hspace*{.3in}(i) $u_{2k-1}(E\oplus F)= u_{2k-1}(E) + u_{2k-1}(F)$;\\
\hspace*{.3in}(ii) $u_{2k-1}(E\otimes F)= rk(E) u_{2k-1}(F)+ rk(F) u_{2k-1}(E)$;\\
\hspace*{.3in}(iii) $u_{2k-1}(E^{*})= - u_{2k-1}(E)$. In particular, if $E$ admits an invariant metric, then $u_{2k-1}(E)$ vanish;\\
\hspace*{.3in}(iv) $u_{2k-1}(E)= 0$ if $k$ is even and $E$ is real.\\

\hspace*{.1in}Combining these invariants $u_{2k-1}$ with the obvious invariant given by the rank of a representation, we
obtain a map:
\begin{equation}\label{theUmap}
U: Rep_{\mathbb{C}}(\mathfrak{g}) \rmap \mathbb{Z}\ltimes H^{odd}(\mathfrak{g}; \mathbb{R}) ,
\end{equation}
and the properties above translate into the fact that $U$ is a morphism of $*$-semi-rings.\\

{\it Proof:} This is only a closer look at the definition. For instance $(iv)$ means that the restriction map
\[ H^{*}(gl_n(\mathbb{C}), U(n)) \rmap H^{*}(gl_n(\mathbb{R}), O(n)) \]
kills the universal $u_{4k-1}$ (see also \cite{Crai2}).  \ $\Boxe$\\

{\it Proof of Theorem \ref{vanishing}:} We now  freely use the language of $\mathfrak{g}$-DG algebras;
an exposition can be found e.g. in Chapter $3$ of \cite{KT}. Let $E$ be a vector bundle over $M$.
Recall (see e.g. \cite{RLF}) that a $\mathfrak{g}$-connection on $E$ is a linear map $\Gamma(\mathfrak{g})\times\Gamma(E)\rmap \Gamma(E)$,
$(X, s)\mapsto \nabla_{X}(s)$, satisfying (\ref{unu}) and (\ref{doi}). It is called flat if $\nabla_{[X, Y]}= [\nabla_X, \nabla_Y]$.
Hence a representation of $\mathfrak{g}$ is a vector bundle over $M$ endowed with a flat $\mathfrak{g}$-connection.
Exactly as in the flat case above, we find a $1-1$ corespondence between $\mathfrak{g}$-connection and connection
$1$-forms $\omega= \omega_{\nabla}\in C^1(\pi^*\mathfrak{g})\otimes gl_n$ on the $gl_n$-DG algebra $C^*(\pi^*\mathfrak{g})$. The Chern-Weil construction
gives us a characteristic
map $k_{\nabla}: W(gl_n)= C^*(gl_n)\otimes S(gl_{n}^{*}) \rmap C^*(\pi^*\mathfrak{g})$ defined on the Weil complex of $gl_n$.
Passing to $gl_n$-basic elements, it induces a map $Inv(gl_n)\rmap C^*(\mathfrak{g})$ defined on the algebra
of invariant symmetric polynomials on $gl_n$. The map induced in cohomology $k_{\nabla}: Inv(gl_n) \rmap  H^*(\mathfrak{g})$ 
will be independent of the connection $\nabla$. In particular, one can use a $\mathfrak{g}$- connection which is induced
by an usual connection on the vector bundle $E$, and we see that $k_{\nabla}$ is just the composition
of the map $\rho_{*}: H^{*}(M)\rmap H^{*}(\mathfrak{g})$ with the usual characteristic map $k_{E}: Inv(gl_n)\rmap H^{*}(M)$
of $E$ (and which defines the Chern classes of $E$). On the other hand, if $\nabla$ is a flat $\mathfrak{g}$-connection on
$E$, then $k_{\nabla}$ kills the symmetric part of $W(gl_n)= C^*(gl_n)\otimes S(gl_{n}^{*})$, hence the map in cohomology
is trivial.  \ $\Boxe$\\

\begin{exam}\label{example3}\emph{ When $\mathfrak{g}= TM$ we obtain the usual characteristic classes (and their construction)
of flat vector bundle (see e.g. \cite{KT}). \\
\hspace*{.3in}When $\mathfrak{g}= \F$ is a foliation of a manifold $M$, any foliated bundle $E\in Rep(\F)$
defines a vector bundle $E_{L}= E|_{L}$ endowed with a flat connection, for each leaf $L$ of $\F$. Similarly, 
the foliated cohomology $H^*(\F)$ can be viewed as a glueing of the
cohomology groups $\{H^*(L)\}_{L}$, and $u_{2k-1}(E)$ can be viewed as
a glueing of the usual characteristic classes $u_{2k-1}(E_{L})\in H^{2k-1}(E_{L})$ of the flat vector bundles $E_{L}$. 
When applied to the normal bundle $\nu= TM/\F$ of the foliation, we obtain certain characteristic
classes $u_{2k-1}(\nu)\in H^{2k-1}(\F)$, $1\leq i\leq q$ where $q$ is the codimension of $\F$. Note the intimate relation with the secondary characteristic classes found by Bott \cite{Bo}, living in $H^*(M)$. What happens is that the image of those 
classes in $H^*(\F)$ vanish (because the Chern classes vanish), and $u_{2k-1}(\nu)$ are precisely the new relevant classes that live in $H^*(\F)$.}
\end{exam}

\begin{exam}\label{fernandes} {\bf (intrinsic characteristic classes):}\emph{ Similar characteristic classes $u_{2k-1}(\mathfrak{g})$, which depend on  $\mathfrak{g}$ only (and not on an auxiliary representation) have been defined by R.L. Fernandes \cite{RLF}. To describe the relation to our classes, we
first assume that $\mathfrak{g}$ is regular, i.e. the image $\F$ of the anchor map has constant rank. Let $\nu$ be the normal bundle of $\F$,
and let $K$ be the kernel of the anchor map. With the Bott connections (see Examples \ref{Bottex}), $K$ and $\nu$ are representations of $\mathfrak{g}$, and 
one can show \cite{Crai2} that
\[ u_{2k-1}(\mathfrak{g})= u_{2k-1}(K)- u_{2k-1}(\nu) .\]
This suggests that Fernandes' $u_{2k-1}(\mathfrak{g})$ can be viewed as the 
characteristic class of the ``formal difference'' $K- \nu$. Since there are exact sequences of vector bundles
$0\rmap K\rmap \mathfrak{g} \rmap \F\rmap 0$ and $0\rmap \F\rmap TM\rmap \nu\rmap 0$, the previous difference bundle equals
to $\mathfrak{g}- TM$ (view these in the $K$-theory of $M$). Hence the classes $u_{2k-1}(\mathfrak{g})$ are the new secondary 
classes which arise from the following vanishing result (implied by Theorem \ref{vanishing}): }
\begin{lem} For any regular Lie algebroid $\mathfrak{g}$ over $M$, one has $Ch(\mathfrak{g}- TM)= 0$ in $H^{*}(\mathfrak{g})$.
\end{lem}
\emph{This can be viewed as an analogue of Bott vanishing theorem \cite{Bo} for characteristic classes of foliations. Using a
nice adaptation of Bott's methods, Fernandes proves this result (and constructs the resulting secondary classes $u_{2k-1}(\mathfrak{g})$)
without any regularity assumption. Although we cannot extend our interpretation of $\mathfrak{g}- TM$ as a representation of $\mathfrak{g}$
from the regular to the non-regular case, this formal difference is always a ``representation up to homotopy''
of $\mathfrak{g}$ in the sense of \cite{ELW} (called there ``the adjoint representation''). Moreover, we can extend our 
characteristic classes from representations to representations up to homotopy, and the conclusion is that $u_{2k-1}(\mathfrak{g})$ of \cite{RLF} are 
always the characteristic classes of the adjoint representation.
For details see \cite{Crai1, Crai2}.}
\end{exam}

\subsection{Relation with the Van Est map}
\label{reltoVE}

\hspace*{.3in}In this subsection we shortly discuss the relation between the characteristic
classes previously introduced and the differentiable cohomology: we show that
the characteristic map 
\[ U: Rep(\mathfrak{g}) \rmap \mathbb{Z}\ltimes H^{odd}(\mathfrak{g}; \mathbb{R})\]
constructed in the previous section naturally factors, via the Van Est map, through the differentiable cohomology.
More precisely, making use both of our Morita invariance and Van Est isomorphism we show:

\begin{st}\label{RwVE} If $E$ is a representation of a Lie groupoid $\G$, and $\mathfrak{g}$ is the
Lie algebroid of $\G$, then the characteristic classes $u_{2k-1}(E)\in H^{2k-1}(\mathfrak{g})$
lie in the image of the Van Est map $\Phi: H^{*}_{d}(\G)\rmap H^{*}(\mathfrak{g})$.
\end{st}

\begin{exam}\emph{The $\alpha$-simply connected integration of the tangent Lie algebroid $TM$ of a manifold 
is the fundamental groupoid of $M$, which is Morita equivalent to the fundamental group $\pi(M)$ of $M$ (cf. Example \ref{example1}). 
Hence, by the Morita invariance of differentiable cohomology, we see that our
Theorem \ref{RwVE} becomes the well-known result that characteristic classes of flat vector bundles come from the
cohomology $H^*(\pi(M))$ of the discrete group $\pi(M)$ (see also Example \ref{example3}).\\
\hspace*{.3in}Given a foliation $(M, \F)$, since the normal bundle of $\F$ is endowed with an action
of the holonomy groupoid $Hol(M, \F)$ (see e.g. \cite{Crath} for more on representations
of the holonomy groupoid), the previous theorem tells us that the classes $u_{2k-1}(\nu)$ (see Example \ref{example3})
come from $H^{*}_{d}(Hol(M,\F))$. In particular, they vanish if the leaf space is an orbifold. More generally, using 
Proposition \ref{properth}, we see that if
a Lie algebroid admits a proper $\alpha$-simply connected integration, then all the characteristic classes of its representations
must vanish.}
\end{exam}

{\it Proof of Theorem \ref{RwVE}:} The idea is quite simple. The frame bundle $P= P(E)$ of an $n$-dimensional representation
$E$ of $\G$ is a left $\G$-space and a principal $GL_n$-bundle, hence can be viewed as a generalized
map $\phi_E: \G\rmap GL_n$ (and this defines a $1-1$ correspondence $Rep_n(\G)\cong Hom_{gen}(\G, GL_n)$).
Hence, by the Morita invariance of differentiable cohomology (Theorem \ref{Moritatheorem}) with trivial coefficients
we obtain a map:
\[ \phi_{E}^{*}: H^{*}_{d}(GL_n) \rmap H^{*}_{d}(\G) .\]
On the other hand, by the classical Van Est isomorphism for $GL_n$, and by the computation of $H^*(gl_n, U(n))$
used also in the previous section, $H^{*}_{d}(GL_n)$ is isomorphic to the exterior algebra on generators 
$u_{2k-1}$, $1\leq k\leq n$. What happens is that 
\begin{equation}\label{reltoVEch}
\phi_{E}^{*}(u_{2k-1})= u_{2k-1}(E) .
\end{equation}
Here are the details. The main problem is to compare the various double complexes which are involved
in the definition of our objects. All these double complexes are of Van Est type (below we will fix the notations),
with one exception: the double complex arising in the construction of $\phi_{E}^{*}$. To avoid working with
this double complex, we use the following trick. We consider the pull-back $\pi^*\G$ which is a groupoid over
$P$ whose space of arrows consists of triples $(p, g, q)\in P\times \nG{1}\times P$ with $\pi(q)= \beta(g)$,
$\pi(q)= \alpha(g)$ (with the first and last projections as target and source map, respectively, and with
the obvious composition). For later use, note that its Lie algebroid is precisely $\pi^{!}\mathfrak{g}$
(see Examples \ref{pull-back}).
We will make use if the following (true!) morphisms of groupoids
\begin{equation}\label{thesequence} 
\G \stackrel{f_{\pi}}{\lmap} \pi^*\G \stackrel{u_{\pi}}{\rmap} GL_n .
\end{equation}
Here $f_{\pi}$ is the obvious projection, while $u_{\pi}$ associates to an arrow $(p, g, q)$ of $\pi^*\G$
the unique matrix $A$ such that $pA= gq$. The main property of these morphisms
is that $f_{\pi}$ is an essential equivalence (i.e. defines a Morita equivalence; see \ref{basicgroupoids}), and, as
generalized morphisms, $\phi_{E}= u_{\pi}f_{\pi}^{-1}$.\\
We now look at the Van Est maps. For any Lie groupoid $\G$, denote by $C(\G)$ the double complex appearing in the proof of Theorem \ref{theWeXu}
(denoted $C$ there); recall that it connects the differentiable and the algebroid cohomology:
\begin{equation}\label{three} 
C_{d}(\G) \stackrel{\epsilon}{\rmap} C(\G) \stackrel{\eta}{\lmap} C(\mathfrak{g}) 
\end{equation}
Moreover, $\eta$ was a quasi-isomorphism, and the map in cohomology $\epsilon \eta^{-1}$ is
the Van Est map for $\G$. For $\G= GL_n$, both $\epsilon$ and $\eta$ are quasi-isomorphism 
if one passes to $U_n$-basic elements in the last two complexes of (\ref{three}), and this
describes the classical Van Est isomorphism for $GL_n$, $H^{*}_{d}(GL_n)\cong H^{*}(gl_n, U_n)$.
Since (\ref{three}) is natural on $\G$ with respect to morphisms of groupoids, by applying it
to each of the groupoids in (\ref{thesequence}), we obtain a commutative diagram
\[ \xymatrix{
C_{d}^{*}(\G)  \ar[d]\ar[r]^-{\sim} & C_{d}^{*}(\pi^{!}\G) \ar[d] & C_{d}^{*}(GL_n) \ar[d]_-{\epsilon}\ar[l] \\
C(\G)  \ar[r]^-{f_{\pi}^{*}} & C(\pi^{!}\G)  & C(GL_n) \ar[l] \\
C^{*}(\mathfrak{g})\ar[u]^-{\sim} \ar[r]^-{f_{\pi}^{*}} & C^{*}(\pi^{!}\mathfrak{g}) \ar[u]^-{\sim} & C^{*}(gl_n)\ar[l]_-{k_{E}}\ar[u]^-{\sim} } \]
Maps marked with ``$\sim$'' are those which are quasi-isomorphisms, and it is not difficult to see
that $k_{E}$ is precisely the map (\ref{kappaE}) used to construct the characteristic classes of $E$. 
By passing to $U(n)$-basic elements in the complexes which form the small diagram in the bottom right corner,
we obtain a similar diagram with the additional property that the maps denoted above by $f_{\pi}^{*}$ and $\epsilon$ 
become quasi-isomorphisms (cf. Theorem \ref{fibers}). Passing to cohomology, and using the previous remarks about the Van Est maps, about the 
map $\phi_{E}^{*}$, and about the definition of $k_{E}$ in (\ref{thecharactr}), we obtain a commutative diagram
\[ \xymatrix{ 
H^{*}_{d}(\G) \ar[d]_-{\Phi} & H^{*}_{d}(GL_n) \ar[d]_-{\sim}\ar[l]_-{\phi_{E}^{*}} \\
H^{*}(\mathfrak{g}) & H^{*}(gl_n, U_n) \ar[l]_-{k_{E}} } \]
which concludes the proof of the theorem. \ $\Boxe$\\

\section{Applications to Poisson manifolds}
\label{Poiss}

\hspace*{.3in}In this section we discuss some applications of our results to Poisson manifolds. Although
this applications are sort of obvious, our intention is to show the relevance of Lie algebroids and of our results
to those interested on Poisson geometry, but less interested on the general theory of Lie groupoids/algeborids.
As a resum\'e: we derive a new proof (and a slight improvement) of the well-known Dazord-Hector integrability criterion for Poisson manifolds \cite{DaHe},
we clarify the problem of Morita invariance of Poisson cohomology (known in degree one only \cite{GiLu}), we prove the Morita invariance
of the modular class (known under certain conditions only \cite{Gi}), we explain the nature of the modular class for regular
Poisson manifolds, and we argue that the first Poisson cohomology groups and our characteristic classes are obstructions to a representation theory for Poisson manifolds
which is analogous to the representation theory for compact Lie groups.\\

\hspace*{.1in}For an introduction to Poisson geometry we recommend  \cite{Vais}, as well as \cite{McXu, We, XuM}. 
Recall here that a Poisson manifold is a manifold $P$ together with a $2$-tensor $\pi\in \Gamma(\Lambda^2TP)$
with the property that $\{f, g\}= \pi(df, dg)$ defines a Lie bracket on $C^{\infty}(P)$ satisfying the 
Leibniz identity $\{f, gh\}= \{f, g\}h+ g\{f, h\}$. Call representation of $(P, \pi)$ any vector bundle $E$ over $P$ together
with an external (bilinear) bracket $\{\cdot, \cdot\}: C^{\infty}(P)\times \Gamma(E)\rmap \Gamma(E)$ which satisfies the Leibniz identities
$\{fg, s\}= f\{g, s\}+ g\{f, s\}$, $\{f, gs\}= \{f, g\}s+ g\{f, s\}$,
and the Jacobi identity $\{\{f, g\}, s\}= \{ f, \{ g, s\}\}+ \{ g, \{ f, s\}\}$, for all $f, g\in C^{\infty}(P)$ and $s\in \Gamma(E)$.
Denote by $Rep(P, \pi)$ the semi-ring of (isomorphism classes of)
representations of $P$.\\
\hspace*{.3in}Recall also that associated to $(P, \pi)$ there is a Lie algebroid
$(T^*P, \rho, [\cdot, \cdot])$, where $\rho: T^*P\rmap TP$ is the map induced by $\pi$ (i.e. $\rho(df)$ is the vector
field induced by the derivation $\{f, \cdot\}$), and $[\alpha, \beta]= L_{\rho(\alpha)}\beta- L_{\rho(\beta)}\alpha- \pi(\alpha, \beta)$.
The relevance here of the Lie algebroid $T^*P$ is that the resulting cohomology is isomorphic to the Poisson cohomology of $(P, \pi)$, usually denoted by $H^{*}_{\pi}(P)$ (see \cite{WeXu}), while its representations are precisely the representations of $(P, \pi)$. \\
\hspace*{.3in}One says that $(P, \pi)$ is integrable if there exists a symplectic (Hausdorff)
groupoid over $P$ which induces the given Poisson structure. If $\G$ is the $\alpha$-simply connected integration of
$(P, \pi)$ it follows that $Rep(P, \pi)\cong Rep(\G)$. For more details, as well as for an exposition 
of integrability results, we refer to Chapter $9$ of \cite{Vais}. Using the  fundamental result of Mackenzie-Xu (Theorem $5.2$ 
in \cite{McXu} which states that $P$ is integrable if and only if its Lie algebroid is
integrable by a Hausdorff groupoid), our Theorem \ref{integration} immediately implies the
following slight improvement of the Hector-Dazord integrability result \cite{DaHe} for
regular Poisson manifolds with totally aspherical symplectic foliation.

\begin{cor}\label{corHD} Let $P$ be a regular Poisson manifold, and let $\F$ be the associated
foliation. If:\\
\hspace*{.3in}$(i)$ $\F$ has no (non-trivial) vanishing cycles,\\
\hspace*{.3in}$(ii)$ for any leaf $L$ of $\F$, $\pi_2(L)$ contains only elements of finite order,\\
then $P$ is integrable.
\end{cor}

{\it Proof:} The first condition is equivalent to the Hausdorffness of the monodromy
groupoid $\G$ of the foliation $\F$, while the second one is equivalent
to the vanishing of the second cohomology groups of its $\alpha$-fibers (i.e. of the universal covers
of the leaves of $\F$). Since the kernel of the anchor map $\rho: T^*P\rmap \F\subset TP$
is always abelian \cite{Vais}, we can apply Theorem \ref{integration}.  \ $\Boxe$\\

\hspace*{.1in}Let us now explain how our results clarify the Morita invariance of Poisson cohomology
(known to hold in degree one only) in all degrees. Recall \cite{XuM} that a Morita equivalence 
between two Poisson manifolds $(P_1, \pi_1)$ and $(P_2, \pi_2)$ is a complete full dual pair \cite{We} 
$P_1 \stackrel{\sigma_1}{\lmap} X\stackrel{\sigma_2}{\rmap} P_2$ with connected and simply-connected fibers. 
Morita equivalence only makes sense on the class of
integrable Poisson manifolds, and there it does define an equivalence relation (cf. the Remark on pp. 496, and Corollary 3.1 of \cite{XuM}).

\begin{cor}\label{abcMP} Let $(P_1, \pi_1)$ and $(P_2, \pi_2)$ be two Poisson manifolds. Any Morita equivalence $P_1 \lmap X\rmap P_2$ with homologically $n$-connected fibers induces isomorphisms
\begin{equation}\label{isMP} 
H^{k}_{\pi_1}(P_1) \cong H^{k}_{\pi_2}(P_2) 
\end{equation}
in all degrees $k\leq n$.
\end{cor}

\hspace*{.1in}Since Morita equivalent Poisson manifolds have Morita equivalent $\alpha$-simply connected integrating groupoids (Theorem $3.2$ in \cite{XuM}),
whose $\alpha$-fibers are precisely the fibers of $\sigma_1$ and $\sigma_2$, this result immediately follows from our Morita
invariance of differentiable cohomology (Theorem \ref{Moritatheorem}), and Theorem \ref{theWeXu}. 
Alternatively, one remarks that the pull-back groupoids (see Example \ref{pull-back})
$\sigma_{1}^{!}T^*P_1$ and $\sigma_{2}^{!}T^*P_2$ are isomorphic, and invoke Theorem \ref{fibers}.
Apart from being much more direct, this last argument also show that Poisson cohomology is invariant
under a much weaker notion of Morita equivalence, which does not exclude non-integrable Poisson manifolds (we hope to make
this more clear in some other place). \\
\hspace*{.3in}The first part of our Theorem \ref{fibers} implies that Morita equivalence of Poisson manifolds behaves
as it should:

\begin{cor}\label{corolarpatru}Any Morita equivalence $P_1 \lmap X\rmap P_2$ of Poisson manifolds induces
isomorphisms  
\begin{equation}\label{MorRep}
Rep(P_1, \pi_1)\cong Rep(P_2, \pi_2)
\end{equation}
\end{cor}

\hspace*{.1in}Let me sketch now an argument which explains the nature of isomorphism between the algebroids $\sigma_{1}^{!}T^*P_1$ and $\sigma_{2}^{!}T^*P_2$, isomorphism which has been essential to our conclusions. 
Let $\F(\sigma_1)$ and $\F(\sigma_2)$ be the foliations on $X$ induced by the fibers
of $\sigma_1$ and $\sigma_2$, respectively. Since $\F(\sigma_2)\cong \pi_{1}^{*}T^*P_1$, there is a natural action
of $\F(\sigma_2)$ on $\F(\sigma_1)$ (see Example \ref{exam2}), and, similarly, an action of  $\F(\sigma_1)$ on $\F(\sigma_2)$.
In the terminology of \cite{McKbi}, $(\F(\sigma_1), \F(\sigma_2))$ form a matched pair, and we obtain a Lie algebroid structure
on $\F(\sigma_1)\oplus \F(\sigma_2)$, which we denote by $\F(\sigma_1)\dcross \F(\sigma_2)$. It is quite straightforward to see
now that this Lie algebroid is isomorphic to our $\sigma_{i}^{!}T^*P_i$'s, $i\in \{1, 2\}$. Hence, in the previous two corollaries,
the cohomology/representation of our Poisson manifolds are isomorphic also to the cohomology/representation of $\F(\sigma_1)\dcross \F(\sigma_2)$.\\

\hspace*{.1in}We now turn to the Morita invariance of the modular class of a Poisson manifold, which is known under
certain conditions only (cf. section $4$ in \cite{Gi}). 

\begin{cor}\label{corolar3}Any Morita equivalence $P_1\lmap X \rmap P_2$ between Poisson manifolds induces
an isomorphisms $H^{1}_{\pi_1}(P_1) \cong H^{1}_{\pi_2}(P_2)$ which maps the modular class of $P_1$
into the modular class of $P_2$.
\end{cor}

{\it Proof:} As mentioned in Examples \ref{example2}, the modular class of a Poisson manifold $P$ is just the first characteristic class $u_1(L)$
of a canonical one dimensional representation $L$ (denoted $Q_{A}$ in \cite{ELW}) of the Lie algebroid $T^*P$. 
 Let $G_1$and $G_2$ be symplectic integrations of $P_1$, and $P_2$, respectively, and let
$\phi_{X}: G_1\cong G_{2}$ be the Morita equivalence of groupoids induced by $X$ \cite{XuM}. Then the canonical representations
of $P_1$ and $P_2$ are related by $L_1= \phi_{X}^{*}L_2$, hence it suffices to use Theorem \ref{RwVE} (more precisely
the formula (\ref{reltoVEch}) appearing in the proof), and the naturality ensured by Theorem \ref{Moritatheorem}.  
As indicated in the previous discussions, one can also give a direct proof at the infinitesimal (algebroid) level. \ $\Boxe$\\

\hspace*{.1in}Clearly, the same argument applies also to the higher characteristic classes.
By the constructions of this section we obtain, for any $E\in Rep(P, \pi)$, certain classes $u_{2k-1}(E)\in H^{*}_{\pi}(P)$, which fit into
a homomorphism
\[ U: Rep(P, \pi) \rmap \mathbb{Z}\ltimes H^{odd}_{\pi}(P) .\]

\begin{cor}\label{corolar4}Let $P_1 \lmap X\rmap P_2$ be a Morita equivalence of Poisson manifolds
with homologically $n$-connected fibers, and let $2k-1\leq n$. For any two representations $E_1$, $E_2$ which correspond to each other
in the Morita isomorphism $Rep(P_1, \pi_1)\cong Rep(P_2, \pi_2)$, the isomorphism  
\[ H^{2k-1}_{\pi_1}(P_1) \tilde{\rmap} H^{2k-1}_{\pi_2}(P_2) \]
(cf. Corollary \ref{abcMP}) maps $u_{2k-1}(E_1)$ into $u_{2k-1}(E_2)$.
\end{cor}

{\it Proof:} The Morita invariance of representations follows from the corresponding result for groupoids
(see subsection \ref{DcMi}), and Theorem 3.2 of \cite{XuM} again.  \ $\Boxe$\\

\hspace*{.1in}Now, a word about the modular class of regular Poisson manifolds $(P, \pi)$. Let $\F$ be the
symplectic foliation of $P$, and let 
\begin{equation}\label{anchorincoh} 
\rho^{*}: H^{*}(\F)\rmap H^{*}_{\pi}(P) 
\end{equation}
be the map induced in cohomology by the anchor. The following result shows that mild conditions on the geometry of the symplectic foliation
(e.g. vanishing of $H^{1}(\F)$) implies the vanishing of the modular class of $P$. Note the formula (\ref{formodP}) in the proof below, which expresses the modular class in terms of our characteristic classes.

\begin{cor}\label{modinfol}The modular class of a regular Poisson manifold $(P, \pi)$ comes from the foliated cohomology of the symplectic
foliation $\F$ of $\pi$, i.e. is in the image of (\ref{anchorincoh}).
\end{cor}

{\it Proof:} Denote by $\mathfrak{g}$ the Lie algebroid of $P$, and let $K$, $\nu$ and $\F$ be as in Examples \ref{example3}. 
By the arguments given there, $K$ and $\nu$ are representations of $\mathfrak{g}$. Up to multiplication to a constant, the modular class of $P$
equals to
\begin{equation}\label{formodP}
mod(P)= u_1(K)- u_1(\nu)\ \in H^{1}_{\pi}(P) \ .
\end{equation} 
Indeed, choosing a linear splitting of the short exact sequences of vector bundles relating $K$, $\nu$, $\F$, $\mathfrak{g}$, and $TM$, we
obtain isomorphisms $TP\cong \F\oplus \nu$, $\mathfrak{g}\cong K\oplus \F$. The induced isomorphisms
$\Lambda^{top}TP\cong \Lambda^{top}\F\otimes \Lambda^{top}\nu$, and $\Lambda^{top}\mathfrak{g}\cong \Lambda^{top}K\otimes \Lambda^{top}F$
are natural, i.e. do not depend on the choice of the linear splittings. It follows that, as vector bundles,
$Q_{\mathfrak{g}}= \Lambda^{top}K \otimes \Lambda^{top}\nu^{*}$. It is not difficult to see
that the action of $\mathfrak{g}$ on $Q_{\mathfrak{g}}$ \cite{ELW} comes from the canonical actions of $\mathfrak{g}$
on $K$ and on $\nu$ (mentioned above), hence $u_{1}(Q_{\mathfrak{g}})= u_{1}(K)- u_{1}(\nu)$. It suffices to remark that $K$ and $\nu \in Rep(P, \pi)$ are actually representations of $\F$ (of Bott type cf. our Examples \ref{Bottex}).  \ $\Boxe$\\

\hspace*{.1in}The previous argument obviously applies to the higher intrinsic characteristic classes of $P$ as well
(but not to those of general representations). Regarding these classes, we mention here the following immediate
consequence (which is known in the case of the modular class)

\begin{cor} If a representation $E\in Rep(P, \pi)$ admits a Poisson-invariant metric, then all its characteristic classes
vanish.
\end{cor}

\hspace*{.1in}The Poisson invariance of a metric $h$ on $E$ is just:
\[ \{ f, h(s_1, s_2)\}= h(\{ f, s_1\}, s_2)+ h(s_1, \{ f, s_2\}), \ \ \forall f\in C^{\infty}(P), s_1, s_2\in \Gamma(E) .\]
Related to the previous corollary, we mention the following result which shows that the first Poisson cohomology group, together with
our characteristic classes are a (serious, we believe) obstruction to a close relation between the geometry of Poisson manifolds
and the one of compact Lie groups.

\begin{cor}If the $\alpha$-simply connected integration of a Poisson manifold $(P, \pi)$ is a proper groupoid, then: \\
\hspace*{.3in}(i) $H^{1}_{\pi}(P)= 0$;\\
\hspace*{.3in}(ii) The characteristic classes of any representation $E\in Rep(P, \pi)$ vanish;\\
\hspace*{.3in}(iii) Any representation of $(P, \pi)$ admits a Poisson-invariant metric.
\end{cor}

\hspace*{.1in}Clearly $(i)$ and $(ii)$ follow from Proposition \ref{properth}, Theorem \ref{theWeXu} (case $n=1$), and Theorem \ref{RwVE}. Note that $(ii)$ also follows from $(iii)$ and from the general properties of characteristic classes mentioned in the previous section.
The last part $(iii)$ follows by a classical averaging argument, using a cut-off function as in the proof of Proposition
\ref{properth}. \\
\hspace*{.3in}Note also that $(iii)$ implies, by usual the arguments of the representation theory of compact Lie groups,
that any representation can be written as a direct sum of irreducible ones, and rises many questions regarding this analogy.

Marius Crainic,\\
\hspace*{.2in}Utrecht University, Department of Mathematics,\\
\hspace*{.2in}P.O.Box:80.010,3508 TA Utrecht, The Netherlands,\\ 
\hspace*{.2in}e-mail: crainic@math.ruu.nl\\
\hspace*{.2in}home-page: http://www.math.uu.nl/people/crainic/

\end{document}